\documentclass[12pt]{amsart}
\usepackage[utf8]{inputenc}
\usepackage{amsmath}
\usepackage{amssymb}
\usepackage{tikz-cd}
\usepackage{url}
\usepackage{subcaption}
\usepackage{mathtools}

\usepackage{enumitem}
\usepackage{amsthm}
\usepackage{thmtools, thm-restate}
\usepackage{hyperref}
\usepackage[capitalise]{cleveref}
\usepackage{bibnames}

\newtheorem{theorem}{Theorem}[section]

\newtheorem{lemma}[theorem]{Lemma}
\newtheorem{proposition}[theorem]{Proposition}

\theoremstyle{definition}
\newtheorem{definition}[theorem]{Definition}
\newtheorem{example}[theorem]{Example}%
\newtheorem{remark}[theorem]{Remark}%

\newcommand{\Id}{\operatorname{Id}}
\newcommand{\im}{\operatorname{im}}

\newcommand{\bR}{\mathbb{R}}
\newcommand{\B}{\mathbf{B}}
\newcommand{\kk}{\mathbf{k}}

\newcommand{\wa}{{\widetilde{a}}}
\newcommand{\wb}{{\widetilde{b}}}
\newcommand{\wc}{{\widetilde{c}}}
\newcommand{\wde}{{\widetilde{d}}}

\newcommand{\wt}{{\widetilde{t}}}

\newcommand{\veck}{\mathbf{Vec}_\kk}

\newcommand{\ladder}{\textbf{Ladder}}

\newcommand{\dint}{d_I^{\bR}}
\newcommand{\dintn}[1]{d_I^{\bR^{#1}}}

\newcommand{\gint}[2]{{\langle #1^\bullet, #2^\bullet\rangle}}

\usepackage[british]{babel}
\usepackage[useregional]{datetime2}
\DTMlangsetup[en-GB]{showdayofmonth=false}

\tikzset{
    labl/.style={anchor=south, rotate=270, inner sep=.5mm}
}

\begin{document}

\title{A Stable Multipersistence Invariant for Ladder Modules }


\author[Delgado-Garrido]{Marco Delgado-Garrido}
\address{%
IMUS,
Departamento de Matemática Aplicada I\\
Universidad de Sevilla\\
Seville, Spain}
\email{mdelgado10@us.es}

\author[Salbu]{Lars M. Salbu}
\address{%
Department of Informatics\\
University of Bergen\\
Bergen, Norway}
\email{lars.salbu@uib.no}

\author[Torras-Casas]{Álvaro Torras-Casas}
\address{%
Departamento de Matemática Aplicada I\\
Universidad de Sevilla\\
Seville, Spain}
\email{atorras@us.es}

\author[Gonzalez-Diaz]{Rocio Gonzalez-Diaz}
\address{%
Departamento de Matemática Aplicada I\\
Universidad de Sevilla\\
Seville, Spain}
\email{rogodi@us.es}



\date{\today}


\begin{abstract}
We define an additive functor $S$ that assigns a three-parameter persistence module to every morphism between one-parameter persistence modules. We show that if the morphism is bar-to-bar, then its image under $S$ is rectangle decomposable. We prove that $S$ is stable with respect to the interleaving distance and study further properties of the functor related to rectangle decompositions.

\end{abstract}
\maketitle


\section{Introduction}\label{sec1}

Topological Data Analysis (TDA) has been one of the main active research lines in computational topology during the last decade. Persistent homology \cite{carlsson2020persistent} is the most important tool in TDA and has led to broad applications in experimental sciences such as neuroscience and spatial biology, as well as developments in pure mathematics in fields like symplectic geometry and analysis {\cite{polterovich}. In particular, the study of the algebraic structures underlying persistent homology called persistence modules is still a fruitful area of research.

Persistence modules are functors from a poset viewed as a category to the category of vector spaces over some fixed field.

The cornerstone of the development of the algebraic theory of persistence modules is the decomposition theorem \cite{Botnan2018DecompositionOP}, which, under mild conditions, allows us to decompose persistence modules over a totally ordered set as the direct sum of interval modules in an essentially unique way. The unique multiset of intervals that determines the decomposition of a persistence module is called the barcode.

One of the most prominent research directions in the theory of persistence modules is multiparameter persistence. The motivation behind multiparameter persistence is to develop robust methods to analyze the shape of the data. In this setup, under mild conditions,
we can still apply decomposition results; 
however, the indecomposables are not intervals in general and the representation theory of such modules is of wild type. 
Nonetheless, there are many ways to overcome these theoretical obstacles and to define stable invariants in order to replace the notion of barcode. For an extended presentation of multipersistence we refer to~\cite{botnan2023introduction}.

In this context, our work studies the structure of morphisms between persistence modules, also known as ladder modules. These kinds of modules cannot be decomposed as the direct sum of interval modules in general. There have been a number of invariants that capture much information about these modules, such as the 6 pack \cite{dimontesano2024chromatictopologicaldataanalysis} and the additive partial matchings induced by persistence maps \cite{GonSor25}. However, the 6-pack does not give us direct information about the decomposability of the ladder module and the matching is not stable in general.

We present a new invariant for ladder modules, the functor $S$, which assigns to each ladder module a 3 parameter persistence module with a simple combinatorial structure, in such a way that, if the original ladder module $\alpha$ decomposes as the direct sum of persistence morphisms between interval modules, then $S(\alpha)$ is rectangle decomposable. The advantage of rectangle decomposable modules is that their structure is very simple and they have very good metric properties \cite{Batanlandi}.

We also show that the functor $S$ is an additive invariant and is stable with respect to the interleaving distance, providing a direct connection between ladder modules and multiparameter persistence and the relation between bar-to-bar persistence maps and rectangle decomposable modules.

\subsection{Structure of the paper}
The paper is structured as follows:

Section two provides the background and basic results about persistence modules needed throughout the paper. In Section 3 we define the functors $\im$ and $\ker$ and prove their main properties. Sections 4 and 5 are devoted to define the main invariant called $S$ and provide its main properties. In the last couple of sections we presentthe main results of the paper concerning the interaction between $S$ and the decomposition of $\alpha$ in morphisms between interval modules. Finally, Appendix A and B are devoted to present in full detail the technical proofs of some of the main results of the paper.

\section{Background}

\subsection{Persistence Modules and Ladder Modules}
In this manuscript we fix a field $\kk$. In addition, given a poset $P$ we view it as a category with objects in $P$ and a unique morphism $p\to q$ whenever $p\leq q$. A \emph{persistence module (over $P$)} is a functor $V:P\to\veck$. We denote by $V_p$ the vector space $V(p)$, and the structure maps $V(p\leq p')$ by $V_{p,p'}:V_{p}\to V_{p'}$. Inequalities between elements are unique, so we have $V_{p,p}=\Id_{V(p)}$ and $V_{p,p''}=V_{p',p''}\circ V_{p,p'}$ for all $p\leq p'\leq p''$. A \emph{(persistence) submodule} $U$ of $V$ is a persistence module such that $U_p\subseteq V_p$ and $U_{p,q}=V_{p,q}|_{U_p}$ for all $p\leq q$. 

A \emph{persistence map} is a natural transformation between a pair of persistence modules. In particular, persistence modules over a poset $P$ form a functor category $\veck^P$. This category is abelian as the target category $\veck$ is an abelian category \cite[8.3.6.(ii) and 8.3.7.(i)]{kashiwara_schapira_2006}.
Therefore, given a persistence map $\alpha\colon V\to U$, we may consider its image and kernel submodules: $\im(\alpha) \subset V$ and $\ker(\alpha)\subset U$.

Persistence maps are also called \emph{ladder modules}, and we can define morphisms between them: For persistence maps $\alpha:V\to U$ and $\alpha':V'\to U'$ over the same poset $P$ we define a \emph{ladder map} $\beta:\alpha\to\alpha'$ as a pair $(\beta_V,\beta_U)$ of persistence maps $\beta_V:V\to V'$ and $\beta_U:U\to U'$ such that $\beta_U\circ\alpha = \alpha'\circ\beta_V$. Ladder modules with ladder maps as morphisms define a category $\ladder$.

\begin{remark}
    There is an equivalent definition of ladder modules as persistence modules over the product poset $P\times \{0,1\}$ where $\{0,1\}$ has the total order $0<1$. 
    Indeed, for functors $W:P\times \{0,1\} \to \veck$ we get persistence maps $W_{(-,0),(-,1)}:W_{(-,0)}\to W_{(-,1)}$. 
    Conversely, for a persistence map $\alpha:V_0\to V_1$ we can define $W:P\times \{0,1\} \to \veck$ by letting $W_{(-,i)}=V_i$ with structure maps $W_{(p,i),(p',i)}=(V_i)_{p,p'}$ for $p\leq p'$ and $W_{(p,0),(p,1)}=\alpha_p$. These are inverse constructions, and ladder maps are natural transformations in this setting. 
\end{remark}

For some results, general persistence modules over $\mathbb{R}$ do not have enough structure. There are different classes of well behaved modules (e.g.\ see \cref{remark: structure theorem}), and we will use the following two properties:
\begin{definition}\label{def:isolated-critical-finite-length} Let $\alpha:V\to U$ be a ladder module over $\mathbb{R}$.
\begin{itemize}
    \item The persistence module $V$ has \emph{isolated critical points} if for each $p\in\mathbb{R}$ there exist $\varepsilon>0$ such that $V$ is constant on $[p-\varepsilon, p)$ and on $(p,p+\varepsilon]$.
    \item The ladder module $\alpha$ has \emph{isolated critical points} if for each $p\in\mathbb{R}$ there exist $\varepsilon>0$ such that $\alpha$ is constant on $[p-\varepsilon, p)$ and on $(p,p+\varepsilon]$.
    \item The persistence module $V$ has \emph{bars of finite length} (or \emph{no infinite bars}) if for all $t\in\mathbb{R}$ there exist $r\leq t\leq s$ where $\im V_{r,t}= 0$ and $\ker V_{t,s}=V_t$.
\end{itemize}
\end{definition}
Note that persistence modules $V$ with isolated critical points and bars of finite length does \underline{not} have to be \emph{pointwise finite dimensional (p.f.d)}, i.e.\ the dimension of $V_t$ might still be infinite. 

\subsection{Interleaving Distance}\label{sec: interleaving distance}

In this section, we define interleaving distance between persistence modules over $\bR^n$, and between ladder modules over $\bR$.

A \emph{shift} of a poset $P$ is a family of order-preserving maps $\Sigma=\{\Sigma_r:P\to P\}_{r\in [0,\infty)}$ where $0\leq s\leq r$ implies $\Sigma_{{s}}(p)\leq \Sigma_r(p)$ for all $p\in P$, also $\Sigma_0=\Id_P$ and $\Sigma_s(\Sigma_r(p))=\Sigma_{s+r}(p)$ for all $s,r\geq 0$. The map $\Sigma_r$ is called an \emph{$r$-shift}.

If $V:P\to\veck$ is a persistence module and $\Sigma$ is a shift of $P$, then we define the persistence module $V[\varepsilon] = V\circ \Sigma_\varepsilon$ for $\varepsilon\geq0$. In particular, we have $V[\varepsilon]_p = V_{\Sigma_\varepsilon(p)}$. For a persistence map $\alpha:V\to U$, we denote by $\alpha[\varepsilon]$ the $\varepsilon$-shifted persistence map $V[\varepsilon]\to U[\varepsilon]$ where $\alpha[\varepsilon]_p=\alpha_{\Sigma_\varepsilon(p)}$. The shifting operation defines a functor $(-)[\varepsilon]:\veck^P \to \veck^P$.

\begin{definition}\label{def: respecting_shifts}
Let $V:P\to\veck$ be a persistence map, and $\Sigma$ a shift of $P$. We have canonical morphisms $\tau_{\varepsilon,\varepsilon'}^{V,\Sigma}:V[\varepsilon]\to V[\varepsilon']$ 
for every $\varepsilon'\geq\varepsilon\geq 0$, by $(\tau_{\varepsilon,\varepsilon'}^{V,\Sigma})_p = V_{\Sigma_\varepsilon(p),\Sigma_{\varepsilon'}(p)}$. This is a persistence morphism since, $V[\varepsilon]_{p,q}=V_{\Sigma_\varepsilon(p),\Sigma_\varepsilon(q)}$ for $p\leq q$. 
\end{definition}

\begin{definition}[$\varepsilon$-interleaving]
For $\varepsilon\geq 0$, an \emph{$\varepsilon$-interleaving} between persistence modules $V$ and $U$ over $P$ with respect to a shift $\Sigma$ of $P$ consists of a pair of natural transformations
\begin{equation*}
    \phi:V\to U[\varepsilon], \qquad \psi: U\to V[\varepsilon]
\end{equation*}
such that $\psi[\varepsilon]_p\circ\phi_p=V_{p,\Sigma_{2\varepsilon}(p)}$ and $\phi[\varepsilon]_p\circ\psi_p=U_{p,\Sigma_{2\varepsilon}(p)}$ for all $p\in P$, i.e such that the composition of the interleaving morphisms is the canonical map:
    \[ \psi[\varepsilon] \circ \phi = \tau_{0,2\varepsilon}^{V,\Sigma}, \quad \quad \phi[\varepsilon] \circ \psi = \tau_{0,2\varepsilon}^{U,\Sigma}\]
\end{definition}
We say that persistence modules $V$ and $U$ are \emph{$\varepsilon$-interleaved} (with respect to $\Sigma$) if there exists an $\varepsilon$-interleaving between them, and their \emph{interleaving distance} is 
\begin{equation*}
    d_I^\Sigma(V,U)=\inf \left\{\varepsilon \geq 0 \,\middle|\, \textrm{$V$ and $U$ are $\varepsilon$-interleaved}\right\}.
\end{equation*}

We now define interleaving distances between persistence modules over $\bR^n$ and between ladder modules over $\mathbb{R}$ by specifying the shift. 

\medskip

\textbf{(1)} Between persistence modules $V:\mathbb{R}\to\veck$ we take $\varepsilon$-shifts $r\mapsto r+\varepsilon$, so $V$ is shifted to $V[\varepsilon]$ where $V[\varepsilon]_t = V_{t+\varepsilon}$ and $V[\varepsilon]_{t,t'}=V_{t+\varepsilon,t'+\varepsilon}$. This gives the standard interleaving distance between one-parameter persistence modules, and we denote it by $d^\bR_I$; cf. Section~4.1. from \cite{ChazalCohen-Steiner2009}.

\medskip
\textbf{(2)} Between persistence modules $\bR^n\to\veck$ we take $\varepsilon$-shifts by maps $(r_1,\dots,r_n)\mapsto(r_1+\varepsilon,\dots,r_n+\varepsilon)$. This is the interleaving distance along the diagonal from multiparameter persistence, denoted by $d^{\bR^n}_I$; cf. Section~2.2. from~\cite{Lesnick2015}.

\medskip
\textbf{(3)} Between \emph{ladder modules} $\mathbb{R}\times\{0,1\}\to\veck$ we take $\varepsilon$-shifts $(r,i)\mapsto(r+\varepsilon,i)$. 
Note that $\varepsilon$-shifting a ladder module $\alpha:V\to U$ gives exactly the shifted persistence map $\alpha[\varepsilon]:V[\varepsilon]\to U[\varepsilon]$ with the shift from \textbf{(1)}. We denote this interleaving distance by $d^L_I$.

\subsection{Intervals and Rectangles}

In $\mathbb{R}$, an \emph{interval} is a non-empty subset $I\subseteq\mathbb{R}$ that contains all points between two endpoints  $a,b$ in $\mathbb{R}\cup\{\pm\infty\}$. We use \emph{decorated real numbers} as endpoints (cf. Section~2.4. from~\cite{Chazal2016}) to indicate whether the endpoints are included or not: We denote intervals by $I=\langle a^\bullet,b^\bullet\rangle$, where each endpoint is decorated with a $\bullet$ that is either $+$ or $-$, with the following convention:
\begin{itemize}
    \item If $a,b\in\mathbb{R}$, then
\begin{equation*}
    \langle a^-,b^-\rangle = [a,b), \quad
    \langle a^+,b^-\rangle = (a,b), \quad
    \langle a^-,b^+\rangle = [a,b], \quad
    \langle a^+,b^+\rangle = (a,b].
\end{equation*}
    \item Decorated infinite endpoints are always open, e.g.:
 \begin{equation*}
    \langle(-\infty)^\bullet,b^+\rangle = (-\infty,b], \quad
        \langle a^+,(+\infty)^\bullet\rangle = (a,\infty).
   \end{equation*} 
\end{itemize}

For dimensions $n\geq 1$, a subset $R\subseteq \bR^n$ is a \emph{rectangle} if $R=I_1\times\dots\times I_n$ where each $I_i\subseteq\bR$ is an interval. Given a rectangle $R\subseteq\mathbb{R}^n$, the \emph{rectangle module} \cite{Batanlandi} $\mathbf{k}_R$ is the persistence module over $\bR^n$ where for $p\leq p'$ we have
\begin{equation*}
    (\mathbf{k}_R)_p=
    \begin{cases}
        \mathbf{k} &\textrm{ if }p\in R\\
        0 &\textrm{ if }p\not\in R
    \end{cases}
    \qquad \textrm{and} \qquad
    (\mathbf{k}_R)_{p,p'}=
    \begin{cases}
        \Id_{\mathbf{k}} &\textrm{ if }\{p,p'\}\subseteq R\\
        0 &\textrm{ if } \{p,p'\}\not\subseteq R.
    \end{cases}
\end{equation*}
For dimension one, rectangles are exactly the intervals, and rectangle modules are called \emph{interval modules}. We use the convention that $\kk_\emptyset\coloneqq 0$.

Given two persistence modules $V$ and $U$, the \emph{direct sum} $V\oplus U$ is defined point-wise as $(V\oplus U)_t=V_t\oplus U_t$ with structure maps $(V\oplus U)_{p,p'}=V_{p,p'} \oplus U_{p,p'}$ whenever $p\leq p'$. This definition of the direct sum agrees with the categorical direct sum in the abelian category of persistence modules (see e.g.\ \cite[Ch.\ III.3]{lang2002} and \cite[Ex.\ 3.3.9, Ex.3.3.vi]{riehl2017category}). 

A persistence module $V:\mathbb{R}^n\to\veck$ is \emph{rectangle decomposable} (see e.g Def. 2.10 \cite{Batanlandi}) if it can be written as the direct sum of rectangle modules, i.e.\ $M\cong \bigoplus_{R} \mathbf{k}_R$. If $n=1$, we say $V$ is \emph{interval decomposable}, (see, e.g. Def. 2.8 \cite{Batanlandi}), namely
\begin{equation}\label{eq: interval_decomposition}
        V\cong\underset{\langle a^\bullet,b^\bullet \rangle\in \B(V) }{\bigoplus}\kk_{\langle a^\bullet,b^\bullet\rangle},
\end{equation}
here the intervals in the multiset $\B(V)$ are the \emph{barcodes} (or \emph{bars}) of $V$.

As mentioned in the introduction, rectangle decomposable modules have several properties that make them convenient to work with. For instance, 
the bottleneck distance between them and there is a combinatorial formula to compute it . For a short summary of the main results concerning closed formulas for the distance between rectangle decomposable modules check (see e.g Def. 2.10 \cite{Batanlandi}).

\section{Functoriality of Image and Kernel}

In this section, we show how images and kernels of structure maps can be thought of as functors.
This study is motivated by analogous operators used in~\cite{CrawleyBoevey2015}.

Let $V:\mathbb{R}\to \veck$ be a persistence module and recall that for each $r\leq t$ in $\mathbb{R}$ we have a linear map $V_{r,t}:V_r\to V_t$. Considering the image of these structure maps, we define the \emph{image module} as the map $\im V:\mathbb{R}^2\to\veck$ acting on objects $(r,t)$ by
\begin{equation*}
    (\im V)_{(r,t)} =
    \begin{cases}
        \im V_{r,t} & \textrm{ if } r\leq t\\
        V_t & \textrm{ if } r>t,
    \end{cases}
\end{equation*}
and morphisms $(r,t)\leq(r',t')$ are sent to the restriction of $V_{t,t'}$. In a similar fashion, the \emph{kernel module} $\ker V:\mathbb{R}^2\to\veck$ maps pairs $(t,s)$ to
\begin{equation*}
    (\ker V)_{(t,s)} =
    \begin{cases}
        \ker V_{t,s} & \textrm{ if } t\leq s\\
        0\subseteq V_t & \textrm{ if } t>s,
    \end{cases}
\end{equation*}
and morphisms $(t,s)\leq(t',s')$ to the restriction of $V_{t,t'}$.

\begin{lemma}\label{lemma: image_functor}
    If $V$ is a persistence module over $\mathbb{R}$, then the image module $\im V$ and kernel module $\ker V$ are persistence modules over $\mathbb{R}^2$. 
\end{lemma}
\begin{proof}
    We need to show that morphisms are well-defined. Preservation of compositions and identities is inherited from V.
    
    \textit{Image:} For $(r,t)\leq (r',t')$ we show that $V_{t,t'}((\im V)_{(r,t)})$ is a subspace of $(\im V)_{(r',t')}$: If 
    $r'>t'$, then $(\im V)_{(r',t')}=V_{t'}$ which is the codomain of $V_{t,t'}$. So, assume $r'\leq t'$ and consider $x\in (\im V)_{(r,t)}\subseteq V_t$. We have two cases:
    \begin{itemize}
        \item If $t\leq r'$, then $V_{t,t'}(x)=
        V_{r',t'} \circ V_{t,r'}(x)$.
        \item If $r'<t$, so $r\leq r'<t\leq t'$, then there is $y\in V_{r}$ with $V_{r,t}(y)=x$. Now, 
        \begin{equation*}
            V_{t,t'}(x)=V_{t,t'}\circ V_{r,t}(y)=V_{t,t'}\circ V_{r',t} \circ V_{r,r'}(y)=V_{r',t'}\circ V_{r,r'}(y).
        \end{equation*}
    \end{itemize}
    In either case, $V_{t,t'}(x)$ is in $(\im V)_{(r',t')}=\im V_{r',t'}$. 
    
    \textit{Kernel:} Let $(t,s)\leq (t',s')$. If $t>s$, then $(\ker V)_{(t,s)}=0$ and morphisms are trivial. So, let $t\leq s$ and consider $x\in\ker V_{t,s}$, namely $V_{t,s}(x)=0$:
    \begin{itemize}
        \item If $t'\leq s'$, then $V_{t',s'}\circ V_{t,t'}(x)= V_{t,s'}(x)=V_{s,s'}\circ V_{t,s}(x)=V_{s,s'}(0)=0$, so $V_{t,t'}(x)$ is in $\ker V_{t',s'}$. 
        \item If $t'>s'$, then $(\ker V)_{(t',s')}=0$ and $V_{t,t'}(x) = V_{s,t'}(x)\circ V_{t,s}(x)=V_{s,t'}(0)=0$.
    \end{itemize}
    In either case, $V_{t,t'}(x)$ is in $(\ker V)_{(t',s')}$. 
\end{proof}
\begin{remark}\label{remark: image_kernel_submodules}
    If we fix a startpoint $r\in \mathbb{R}$, then the module $(\im V)_{(r,-)}:\mathbb{R}\to\veck$ sending $t$ to $(\im V)_{(r,t)}$ is a submodule of $V$, since $(\im V)_{(r,t)}\subseteq V_t$ and the structure maps of $\im V$ are restricted from $V$. Similarly, fixing the endpoint $s\in\mathbb{R}$ gives a submodule $(\ker V)_{(-,s)}\subseteq V$.
\end{remark} 

Next, we show that $\im$ and $\ker$ are natural with respect to the input persistence module $V$.

\begin{proposition}\label{prop: im_and_ker_naturality}
    Persistence maps $\alpha:V\to U$ over $\mathbb{R}$ restrict to persistence maps $\im V\to\im U$ and $\ker V\to \ker U$ over $\mathbb{R}^2$. In particular, $\im$ and $\ker$ are functors $\veck^\mathbb{R}\to \veck^{\mathbb{R}^2}$.
\end{proposition}
\begin{proof}
    \textit{Image:} We 
    define $\alpha_{(r,t)}:(\im V)_{(r,t)}\to(\im U)_{(r,t)}$ to be the restriction of $\alpha_t$.
    If $r>t$, then $(\im V)_{(r,t)}=V_t$ and $(\im U)_{(r,t)}= U_t$ and $\alpha$ is unchanged. So, assume $r\leq t$ and let $x\in \im V_{r,t}$, so $V_{r,t}(y)=x$ for some $y\in V_r$. Now, $U_{r,t}\circ \alpha_r(y)=\alpha_t\circ V_{r,t}(y)=\alpha_t(x)$, so $\alpha_t(x)$ is in $\im U_{r,t}$. Thus, $\alpha_{(r,t)}$ is well-defined. Naturality is inherited from the structure maps.

    \textit{Kernel:} Similarly, define $\alpha_{(t,s)}:(\ker V)_{(t,s)}\to(\ker U)_{(t,s)}$ by restricting $\alpha_t$. 
    If $t>s$, then $(\ker V)_{(t,s)}=(\ker U)_{(t,s)}=0$ and $\alpha$ is trivial. We look at the case $t\leq s$. If $x\in\ker V_{t,s}$, then $U_{t,s}\circ\alpha_t(x)=\alpha_s\circ V_{t,s}(x)=\alpha_s(0)=0$, so $\alpha_t(x)$ is in the kernel of $U_{t,s}$, and $\alpha_{(t,s)}$ is well-defined. Naturality is inherited from the structure maps.
\end{proof}

\begin{remark}[Structure Theorem {\cite[Thm.\ 1.2]{CrawleyBoevey2015}}]\label{remark: structure theorem}
    If a persistence module $V:\mathbb{R}_{\leq} \to \veck$ satisfies the \emph{descending chain condition (d.c.c.) on kernels and images}, namely if for all $t$ and sequences $q_1> q_2>\dots$ in $\mathbb{R}$ both chains 
    \begin{equation*}
        V_t\supseteq (\im V)_{q_1,t} \supseteq (\im V)_{q_2,t} \supseteq \dots \quad \textrm{and}\quad 
        V_t\supseteq (\ker V)_{t,q_1}\supseteq (\ker V)_{t,q_2}\supseteq\dots 
    \end{equation*}
    stabilize, then $V$ is interval decomposable. In particular, this holds if $V$ has both isolated critical points and bars of finite length.
\end{remark}

\section{Intersections and Adding Variables: The Functor \texorpdfstring{$S$}{S}}\label{sec: Intersections_and_adding_variables}

In the previous section we defined the functors $\im$ and $\ker$ which capture the birth and death features of a one-parameter persistence module. In this section we will combine the information of these two functors to define a new invariant that captures non-trivial information about ladder modules.

For a persistence module $V:P\to\veck$ consider two submodules $A,B\subseteq V$. The \emph{intersection} of $A$ and $B$ is the submodule $A\cap B\subseteq V$ where $(A\cap B)_p=A_p\cap B_p$. This is well-defined since if $x$ is in $A_p\cap B_p$, then $V_{p,q}(x) = A_{p,q}(x)= B_{p,q}(x)$ is in $A_q\cap B_q$ for all $p\leq q$.

If $\alpha:V\to U$ is a ladder module over $P$, then the \emph{image of} $\alpha$ is the submodule $\alpha(V)\subseteq U$ where $\alpha(V)_p=\im \alpha_p\subseteq U_p$. This is well-defined, since if $\alpha_p(y)=x$, then $U_{p,q}(x)=U_{p,q}\circ\alpha_p(y)=\alpha_q\circ V_{p,q}(y)$ is in $\alpha_q(V_q)$ for $p\leq q$.

In \cref{remark: image_kernel_submodules}, we noted that fixing start- and endpoints $r$ and $s$ we get submodules $\im V_{r,-}$ and $\ker V_{-,s}$ of $V$. If $V$
admits an interval decomposition $\B(V)$, the submodules $\im V_{r,-}$ (resp. $\ker V_{-,s}$) admit an interval decomposition consisting of}the bars of 
$\B(V)$ 
that start (at or) before $r$ 
(resp. end before $s$).
 Our aim here is to construct a module that has this information for all $r$ and $s$, while at the same time allows us to study start- and endpoints independently. 
We achieve this by extending the image and kernel modules into a common three-dimensional space, by adding \emph{dummy variables}:
given a persistence module $V:\mathbb{R}^m\to \veck$, we can increase its dimension by pre-composing with projections $\pi^{n\to m}_{i_1,\dots, i_m}:\mathbb{R}^n\to\mathbb{R}^m$ sending $(x_1,x_2,\dots, x_n)$ to $(x_{i_1},x_{i_2},\dots,x_{i_m})$ for $1\leq i_j \leq n$ with $i_s < i_t$ for $s < t$.

\begin{example}\label[example]{example: extending_V}
    For a persistence module $V:\mathbb{R}\to\veck$, we add a dummy variable by pre-composing with the projection $\pi \coloneq \pi^{2\to1}_2:\mathbb{R}^2\to\mathbb{R}$ sending $(t_1,t_2)$ to $t_2$. Now, $V\circ \pi:\mathbb{R}^2\to\veck$ is a persistence module where $(V\circ\pi)_{(r,t)}= V_t$ and $(V\circ\pi)_{(r,t),(r',t')}=V_{t,t'}$. The image module $\im V$ is a submodule of $V\circ \pi$, since $(\im V)_{(r,t)}\subseteq V_t=(V\circ\pi)_{(r,t)}$ with structure maps 
    \begin{equation*}
        (\im V)_{(r,t),(r',t')}=V_{t,t'}|_{(\im V)_{(r,t)}} = (V\circ\pi)_{(r,t),(r',t')}|_{(\im V)_{(r,t)}}.
    \end{equation*}
\end{example}

We now define the main construction that we are exploring in this paper. 

\begin{definition}    
Given a ladder module $\alpha:V\to U$ over $\mathbb{R}$, we define the persistent module $S(\alpha):\mathbb{R}^3\to\veck$ where
\begin{equation}\label{eq: im_cap_ker_definition}
    S(\alpha)_{(r,t,s)} \coloneqq (\im U)_{(r,t)}\cap \alpha_t((\ker V)_{(t,s)}), 
\end{equation}
and the structure maps are given by restricting $U_{t,t'}$. In particular, $S(\alpha)$ is the intersection of  $(\im U)\circ\pi^{3\to2}_{1,2}$ and $\alpha(\ker V)\circ \pi^{3\to 2}_{2,3}$ as submodules of $U\circ \pi^{3\to1}_2$.
\end{definition}
\begin{example}
    If $\alpha: \kk_\gint{a}{b}\xrightarrow{1}\kk_\gint{c}{d}$ is a non-trivial map between intervals, then $c^\bullet\leq a^\bullet\leq d^\bullet\leq b^\bullet$. If $s\in \langle b^\bullet,\infty)$ and $t\in\gint{a}{d}$, then
    \begin{itemize}
        \item $S(\alpha)_{(s,-,s)}=\im\alpha=\kk_\gint{a}{d}$,
        \item $S(\alpha)_{(-,t,s)}=(\im U)_{(-,t)} = \kk_{\langle c^\bullet,\infty)}$, and
        \item $S(\alpha)_{(t,t,-)}=\phi_t((\ker V)_{(t,-)})=\kk_{\langle b^\bullet,\infty)}$,
    \end{itemize}
    so $S(\alpha)$ captures information about all endpoints, $a^\bullet,\, b^\bullet,\, c^\bullet$ and $d^\bullet$. We will see in \cref{sec: bar-to-bar} that this generalizes to more complicated maps.
\end{example}

We continue to show that $S$ is in fact a functor from the category of ladder modules to the category of persistence modules over $\mathbb{R}^3$.

\begin{proposition}\label{prop: S_functorial}
    Consider ladder modules $\alpha:V\to U$ and $\alpha':V'\to U'$ over $\mathbb{R}$. If $\beta=(\beta_V,\beta_U):\alpha\to\alpha'$ is a ladder map, then $S(\beta):S(\alpha)\to S(\alpha')$ 
    given by $(\beta_U \circ \pi^{3\to1}_2)\vert_{ S(\alpha)}$ is a persistence map over $\mathbb{R}^3$.
\end{proposition}
\begin{proof}
    We show that $S(\beta)$ is well-defined; i.e. that its image lies in $S(\alpha')$. Since $\beta_U$ is a persistence map, we get by \cref{prop: im_and_ker_naturality} that $\beta_U(\im U)\subseteq \im U'$. Furthermore, 
    as $\beta$ is a ladder map, it holds that 
    \begin{equation*}
        (\beta_U)_t\circ\alpha_t((\ker V)_{(t,s)}) = \alpha'_t\circ (\beta_V)_t((\ker V)_{(t,s)}),
    \end{equation*}
    and, by \cref{prop: im_and_ker_naturality}, that $(\beta_V)_t((\ker V)_{(t,s)})\subseteq (\ker V')_{(t,s)}$. In particular, 
    if $y$ is in $\alpha_t((\ker V)_{(t,s)})$, then $(\beta_U)_t(y)$ is in $\alpha'_t((\ker V')_{(t,s)})$. 
    On the other hand, 
    if $z$ is in $(\im U)_{(r,t)}\cap \alpha_t((\ker V)_{(t,s)})$, then $(\beta_U)_t(z)$ is in $(\im U')_{(r,t)}\cap \alpha'_t((\ker V')_{(t,s)})$. Naturality is inherited from $\beta_U$. 
\end{proof}

If $\alpha:V\to U$ is a persistence map, then $\Id_\alpha=(\Id_V,\Id_U)$ is the identity morphism and $S(\Id_\alpha) = \Id_{S(\alpha)}$. Furthermore, if $\alpha:V\to U$, $\alpha':V'\to U'$ and $\alpha'':V''\to U''$ are ladder modules with morphisms $\beta=(\beta_V,\beta_U):\alpha\to\alpha'$ and $\beta'=(\beta_{V'},\beta_{U'}):\alpha'\to\alpha''$, then for $y\in S(\alpha)_{(a,t,b)}$ we have
\begin{equation*}
    S(\beta'\circ\beta)_t(y) = (\beta_{U'}\circ\beta_U)_t(y)=(\beta_{U'})_t\circ(\beta_U)_t(y)=S(\beta')_t\circ S(\beta)_t(y).
\end{equation*}
Thus, $S:\ladder\to \veck^{\mathbb{R}^3}$ defines a functor from ladder modules to the category of persistence modules over $\mathbb{R}^3$.

\section{Stability Results for \texorpdfstring{$S$}{S}, \texorpdfstring{$\im$}{im} and \texorpdfstring{$\ker$}{ker}}
\label{sec: Stability}

In this section, we show that the functors $S$,  $\im$ and $\ker$ are stable with respect to the interleaving distance. We will introduce a technical result which simplifies the proofs of the stability properties for $S$, $\im $ and $\ker$:

\begin{definition}
Let $\Sigma^1$ and $\Sigma^2$ be shifts of $P_1$ and $P_2$ respectively. A functor $F:\veck^{P_1}\to\veck^{P_2}$ \emph{respects the shifts $\Sigma^1$ and $\Sigma^2$} if there is a natural transformation $N^\varepsilon:F((-)[\varepsilon])\to F(-)[\varepsilon]$ for each $\varepsilon \geq 0$, i.e.\
\begin{equation*}
    \begin{tikzcd}
        {F(V[\varepsilon])}\ar[d,"{F(\alpha[\varepsilon])}"]\ar[r,"N^{\varepsilon}_V"] 
        & {(F(V))[\varepsilon]}\ar[d,"{F(\alpha)[\varepsilon]}"]\\
        {F(U[\varepsilon])}\ar[r,"N^{\varepsilon}_U"] & {(F(U))[\varepsilon]},
    \end{tikzcd}
\end{equation*}
compatible in the sense that  for all $V\in \veck^{P_1}$ we have
\begin{itemize}
    \item[(1)]\label{def1} $N^0_V=id_{F(V)}$, and
   
    \item[(2)]\label{def3} $N^{\epsilon + \delta}_V=N_V^\epsilon[\delta]\circ N_{V[\epsilon]}^\delta$ for all $\varepsilon,\delta\geq 0$, i.e.\ commuting diagrams
\[\begin{tikzcd}
	{F(V[\varepsilon+\delta])} & {F(V[\varepsilon])[\delta]} & {} \\
	& {F(V)[\varepsilon+\delta]}
	\arrow["{N^{\delta}_{V[\varepsilon]}}", from=1-1, to=1-2]
	\arrow["{N_V^{\varepsilon+\delta}}"', from=1-1, to=2-2]
	\arrow["{N^{\varepsilon}_V[\delta]}", from=1-2, to=2-2].
\end{tikzcd}\]
\item[(3)]  The collection $\{N^\varepsilon_V\}_{\varepsilon\geq 0}$ is compatible with respect to the canonical morphisms $F(\tau_{\varepsilon,\varepsilon'}^{V,\Sigma^1})$ and $\tau_{\varepsilon,\varepsilon'}^{F(V),\Sigma^2}$ for $\varepsilon\leq \varepsilon'$:
\[\begin{tikzcd}
	{F(V[\varepsilon])} & {F(V)[\varepsilon]} \\
	{F(V[\varepsilon'])} & {F(V)[\varepsilon']}
	\arrow["{N^{\varepsilon}_V}", from=1-1, to=1-2]
	\arrow["F(\tau_{\varepsilon,\varepsilon'}^{V,\Sigma^1})"',from=1-1, to=2-1]
	\arrow["\tau_{\varepsilon,\varepsilon'}^{F(V),\Sigma^2}",from=1-2, to=2-2]
	\arrow["{N_V^{\varepsilon'}}"', from=2-1, to=2-2]
    \end{tikzcd}\]
    \end{itemize}   
\end{definition}

\begin{theorem}[Stability of shift respecting functors \label{thm: shift_respecting_stable}]
    Let $\Sigma^1$ and $\Sigma^2$ be shifts of $P_1$ and $P_2$ respectively, and let $V$ and $U$ be persistence modules over $P_1$. If $F:\veck^{P_1}\to\veck^{P_2}$ respects $\Sigma^1$ and $\Sigma^2$, then 
    \begin{equation*}
        d^{\Sigma^1}_I(V,U)\geq d^{\Sigma^2}_I(F(V),F(U)).
    \end{equation*}
\end{theorem}
\begin{proof}
    Proven in Appendix \ref{sec: AppendixA}
\end{proof}
\begin{lemma}\label{lemma: S_im_ker_shift_resp}
    The functors $S$, $\im$ and $\ker$ are shift respecting for the shifts defined in \ref{sec: interleaving distance}.
\end{lemma}
\begin{proof}
    Proven in Appendix \ref{sec: AppendixA}.
\end{proof}

Now we can prove a stability result for the functors $\im$ and $\ker$. In fact, under mild conditions, the functors $\im$ and $\ker$ preserve the interleaving distance:

\begin{restatable}[Isometry for image and kernel\label{thm: isom_im_ker}]{theorem}{isoimker}
For persistence modules $V$ and $U$ over $\mathbb{R}$, we have
    \begin{align*}
        \dint(V,U)=d_I^{\mathbb{R}^2}(\im V,\im U)\quad \textrm{ and }\quad \dint(V,U)\geq d_I^{\mathbb{R}^2}(\ker V,\ker U).
    \end{align*}  
Furthermore, if $V$ and $U$ have bars of finite length, then
\[\dint(V,U)=d_I^{\mathbb{R}^2}(\ker V,\ker U).\]
\end{restatable}
\begin{proof}
    Proven in \cref{sec: AppendixA}.
\end{proof}

Using the previous results stated in this section we can also prove a stability result for the functor $S$:
\begin{theorem}[Stability of $S$\label{thm: stability of S}]
    Let $\alpha:V\to U$ and $\alpha':V'\to U'$ be maps between persistence modules over $\mathbb{R}$. The functor $S$ is stable with respect to interleaving distance, namely
    \begin{equation*}
        d_I^{\mathbb{R}^3}(S(\alpha),S(\alpha'))\leq d_I^L(\alpha,\alpha').
    \end{equation*}
\end{theorem}
\begin{proof}
Applying Lemma \ref{lemma: S_im_ker_shift_resp} and Theorem \ref{thm: shift_respecting_stable} we get the result.
 
\end{proof}
If the ladder modules $V$ and $V'$ are decomposable as a direct sum of finite intervals, we have the following lower bound for the previous inequality:
\begin{proposition}
    Given a pair of ladder modules $\alpha:V\to U , \alpha':V'\to U'$, if $V$ and $V'$ have bars of finite length, then
    \[  \dint (\alpha(V),\alpha'(V')) \leq d_I^{\mathbb{R}^3}(S(\alpha),S(\alpha')).\]
\end{proposition}
\begin{proof}
Proven in Appendix \ref{sec: AppendixA}

\end{proof}

These stability results raise the question wether or not we actually have an isometry like the one that we have for the functors $\im$ and $\ker$. However, as the following example illustrates, this cannot be the case. In fact, the gap in the inequality from \Cref{thm: stability of S} can be arbitrarily big.

\begin{example}
    Consider the two ladder modules $\alpha_1: \kk_{[0,1)}\xrightarrow{1} \kk_{[0,1)}$ and $\alpha_2: \kk_{[0,1+N)}\xrightarrow{1} \kk_{[0,1)}$ where $N>0$. An interleaving of ladder modules restricts to interleavings on domains and codomains, so in particular
    \begin{equation*}
        d^L_I(\alpha_1,\alpha_2)\geq d^\mathbb{R}_I(\kk_{[0,1)},\kk_{[0,1+N)})=N/2.
    \end{equation*}
    However, since $S(\alpha)_{(r,t,s)}\subseteq U_t$, we have $S(\alpha_1)_{(r,t,s)}=S(\alpha_2)_{(r,t,s)}=0$ whenever $t\not\in [0,1)$. Thus, the structure maps $S(\alpha_i)_{(r,t,s),(r+1,t+1,s+1)}$ are all zero for $i=1,2$. So, the trivial morphisms 
    \begin{equation*}
        \phi: S(\alpha_1)\xrightarrow{0} S(\alpha_2)[1/2]\quad\textrm{ and }\quad \psi:S(\alpha_2)\xrightarrow{0} S(\alpha)[1/2]
    \end{equation*}
    define a $1/2$-interleaving: $(\psi[1/2]\circ\phi)_{(r,t,s)}=0=S(\alpha_i)_{(r,t,s),(r+1,t+1,s+1)}$  and similarly for $\phi[1/2]\circ\psi$. We conclude that $d^{\mathbb{R}^3}_I(S(\alpha_1),S(\alpha_2))\leq 1/2$ and
    \begin{equation*}
        d^L_I(\alpha_1,\alpha_2)-d^{\mathbb{R}^3}_I(S(\alpha_1),S(\alpha_2))\geq (N-1)/2
    \end{equation*}
    which can be made arbitrarily large by increasing $N$.
\end{example}

\section{Additivity of \texorpdfstring{$S$}{S}}

In this section, we show that the functor $S$ is additive, using additivity of the image and kernel functors.

A functor $F:\mathcal{A}\to \mathcal{B}$ between abelian categories is \emph{additive} if it preserves direct sums \cite[Thm.\ 3.11]{freyd1964abelian}, i.e. if
    \[ F(0) \cong 0\in \mathcal{B} \quad \text{and } \quad F(X\oplus Y)\cong F(X)\oplus F(Y).\]
Note that the composition of additive functors is additive.

\begin{proposition}\label{proposition: additivity}
The functors $\im$ and $\ker$ are additive. Specifically, if $V,U:\mathbb{R}\to\veck$ are persistence modules, then
\begin{align*}
    \im\,(V\oplus U)&=(\im V)\oplus (\im U),\\
    \ker\,(V\oplus U)&=(\ker V)\oplus (\ker U).
\end{align*}
\end{proposition}
\begin{proof}
Direct sums are defined point-wise, $(V\oplus U)_t=V_t \oplus U_t$, with structure maps $(V\oplus U)_{t,t'}=V_{t,t'}\oplus U_{t,t'}$, and hence we have:
\begin{align*}
    \im\,(V\oplus U)_{a,t}
    &=(\im V_{a,t})\oplus (\im U_{a,t}),\\
    \ker\,
    (V\oplus U)_{t,b}
    &=(\ker V_{t,b})\oplus (\ker U_{t,b}).
\end{align*}
concluding the result.
\end{proof}

\begin{proposition}\label{prop: S_additive}
    The functor $S$ is additive. Specifically, if $\alpha:V\to U$ and $\alpha':V'\to U'$ are ladder modules over $\mathbb{R}$, then
    \begin{equation*}
        S(\alpha\oplus\alpha')=S(\alpha)\oplus S(\alpha).
    \end{equation*}
    \end{proposition}
\begin{proof}
    Direct sums are defined pointwise, so let $r,t,s\in\mathbb{R}$. Using \cref{proposition: additivity} we have
    \begin{align*}
        S(\alpha\oplus& \alpha')_{(r,t,s)} =\\
        &=(\im\, (U\oplus U'))_{(r,t)}\cap (\alpha\oplus\alpha')_t((\ker\, (V\oplus V'))_{(t,s)})\\
        &=(\im U)_{(r,t)}\oplus(\im U')_{(r,t)}\cap (\alpha_t((\ker V)_{(t,s)})\oplus\alpha'_t((\ker V')_{(t,s)}))\\
        &=((\im U)_{(r,t)}\cap\alpha_t((\ker V)_{(t,s)}))\oplus((\im U')_{(r,t)}\cap\alpha'_t((\ker V')_{(t,s)}))\\
      &=S(\alpha)_{(r,t,s)}\oplus S(\alpha')_{(r,t,s)} = (S(\alpha)\oplus S(\alpha'))_{(r,t,s)}
    \end{align*}
    as submodules of $U_t\oplus U'_t$.
\end{proof}

\section{Rectangle Decomposition}\label{sec: Rectangle Decomposition}
In this section we consider how the image and kernel act on intervals. Specifically, we show that intervals are sent to rectangles, and thus the image and kernel modules of general persistence modules are rectangle decomposable.

\begin{lemma}\label{lemma: image_is_rectangle}
    For intervals $\gint{a}{b}$ with $a,b \in \mathbb{R}\cup\{\pm\infty\}$, we have 
    \begin{align*}
        \im \kk_{\gint{a}{b}}&\cong \kk_{\langle a^{\bullet},\infty)\times \gint{a}{b}},\\
        \ker \kk_{\gint{a}{b}}&\cong \kk_{\langle a^{\bullet},b^{\bullet}\rangle\times \langle b^\bullet, \infty)}.
    \end{align*}
\end{lemma}
\begin{proof}
    \textit{Image:} Using \cref{remark: image_kernel_submodules}, for $r,t\in \mathbb{R}$ we have $(\im V)_{(r,t)}\subseteq V_t$. In particular, $\im(\kk_{\gint{a}{b} })_{(r,t)}\subseteq (\kk_{\gint{a}{b}})_t$, which is zero whenever $t\not\in \gint{a}{b}$. If $t\in \gint{a}{b} $, then 
    \begin{equation*}
        (\im \kk_{\gint{a}{b}})_{(r,t)} = 
        \begin{cases}
            \kk &\textrm{if }r\in \langle a^{\bullet},\infty)\\
            0 &\textrm{otherwise}.
        \end{cases}
    \end{equation*}
    The structure maps, when non-zero, are restrictions of the identity map $\Id_\kk:\kk\to \kk$ to $\kk$, thus they are identity maps.

    \textit{Kernel:} Since $(\ker V)_{(t,s)}\subseteq V_t$ for all $t,s\in\mathbb{R}$, it follows that $(\ker \kk_\gint{a}{b})_{(t,s)}= 0$ for $t\not\in \gint{a}{b}$. If $t\in\gint{a}{b}$, then 
    \begin{equation*}
        (\im \kk_{\gint{a}{b}})_{(t,s)} = 
        \begin{cases}
            \kk &\textrm{if }s\in \langle b^{\bullet},\infty)\\
            0 &\textrm{otherwise}.
        \end{cases}
    \end{equation*}
    Structure maps are identities when possible.
\end{proof}
Note in particular that, if $b=\infty$, then $\langle b^\bullet,\infty)=\emptyset$ and so, by our convention, we have $\ker \kk_\gint{a}{b}\cong \kk_\emptyset =0$.

\begin{proposition}\label{corollary: imkerdecomp}
    If $V:\mathbb{R}\to\veck$ is an interval decomposable persistence module with barcode $\B(V)$, then $\im V$ and $\ker V$ are rectangle decomposable with the following decompositions:
     \begin{equation*}
         \im V \cong\bigoplus\limits_{\langle a^{\bullet},b^{\bullet}\rangle \in \B(V)}\kk_{\langle a^{\bullet},\infty )\times\langle a^{\bullet},b^{\bullet}\rangle},
    \end{equation*}
    \begin{equation*}
         \ker V\cong\bigoplus\limits_{\langle a^{\bullet},b^{\bullet}\rangle \in \B(V)}\kk_{\langle a^{\bullet},b^{\bullet}\rangle\times \langle b^\bullet,\infty)}.
    \end{equation*}
\end{proposition}
\begin{proof}
    The module $V$ decomposes into intervals modules like in \cref{eq: interval_decomposition}, and each interval is sent to a rectangle of the given form (\cref{lemma: image_is_rectangle}). The result follows from additivity (\cref{proposition: additivity}).
\end{proof}

\begin{example}[Computation of $\im$ and $\ker$]
Consider the persistence module $V=\kk_{[2,4]}\oplus \kk_{[1,5]}$, let us compute $\im V$ and $\ker V$ using Proposition \ref{corollary: imkerdecomp}:
\[\im V\cong \kk_{[2,\infty)\times[2,4]}\oplus \kk_{[1,\infty)\times [1,5]}
\]
\[
\ker V\cong \kk_{[2,4]\times(4,\infty)}\oplus\kk_{[1,5]\times(5,\infty)}
\]

\begin{figure}[ht]
\centering
\begin{subfigure}{0.48\textwidth}
\centering
\begin{tikzpicture}[scale=0.8,yscale=0.8113]

\draw[->] (-0.2,0) -- (5.6,0) node[below] {$r$};
\draw[->] (0,0) -- (0,5.3) node[left] {$t$};

\fill[blue!50,opacity=.45]
    (2,2) rectangle (5.2,4);

\draw (2,2) -- (2,4);
\draw (2,2) -- (5.2,2);
\draw (2,4) -- (5.2,4);

\fill[orange!70,opacity=.60]
    (1,1) rectangle (5.2,5);

\draw (1,1) -- (1,5);
\draw (1,1) -- (5.2,1);
\draw (1,5) -- (5.2,5);

\draw[->] (5.2,2) -- (5.6,2);
\draw[->] (5.2,4) -- (5.6,4);
\draw[->] (5.2,1) -- (5.6,1);
\draw[->] (5.2,5) -- (5.6,5);

\foreach \x in {1,2}
    \draw (\x,0.08) -- (\x,-0.08);

\foreach \y/\lab in {1/1,2/2,4/4,5/5}
    \draw (0.08,\y) -- (-0.08,\y)
          node[left=4pt] {$\lab$};

\node[below] at (1,0) {$1$};
\node[below] at (2,0) {$2$};

\end{tikzpicture}
\end{subfigure}
\hfill
\begin{subfigure}{0.48\textwidth}
\centering
\begin{tikzpicture}[scale=0.8]

\draw[->] (-0.2,0) -- (5.6,0) node[below] {$t$};
\draw[->] (0,0) -- (0,4.3) node[left] {$s$};

\fill[blue!50,opacity=.45]
    (2,1) rectangle (4,3.8);

\draw[dashed] (2,1) -- (4,1);
\draw (2,1) -- (2,3.8);
\draw (4,1) -- (4,3.8);

\fill[orange!70,opacity=.60]
    (1,2) rectangle (5,3.8);

\draw[dashed] (1,2) -- (5,2);
\draw (1,2) -- (1,3.8);
\draw (5,2) -- (5,3.8);

\draw[->] (2,3.8) -- (2,4.15);
\draw[->] (4,3.8) -- (4,4.15);
\draw[->] (1,3.8) -- (1,4.15);
\draw[->] (5,3.8) -- (5,4.15);

\foreach \x in {1,2,4,5}
    \draw (\x,0.08) -- (\x,-0.08);

\foreach \y/\lab in {1/4,2/5}
    \draw (0.08,\y) -- (-0.08,\y)
          node[left=4pt] {$\lab$};

\node[below] at (1,0) {$1$};
\node[below] at (2,0) {$2$};
\node[below] at (4,0) {$4$};
\node[below] at (5,0) {$5$};

\end{tikzpicture}
\end{subfigure}

\caption{\underline{Left}: Visualization of  $(\im V)_{(r,t)}$.  \underline{Right}: Visualization of $(\ker V)_{(t,s)}$.}
\end{figure}
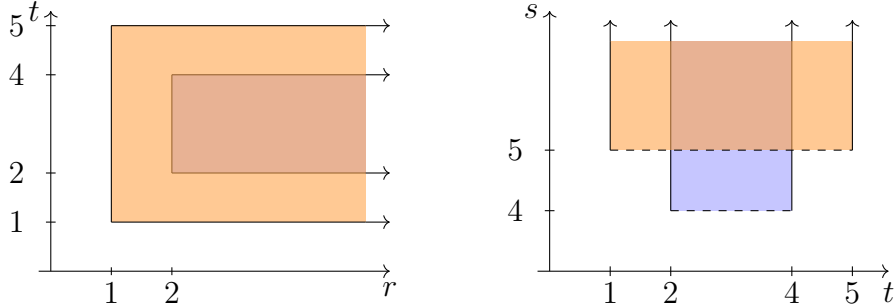
\end{example}

We end this section with an isometry result for the functors $\im$ and $\ker$. For a detailed proof of the following result and the definition of the bottleneck distance for rectangle decomposable modules, we refer the reader to.
\cref{sec: AppendixA}.

\section{Bar-to-Bar}\label{sec: bar-to-bar}
In this section, we show that the functor $S$ defined in \cref{sec: Intersections_and_adding_variables} is rectangle decomposable whenever $\alpha$ can be decomposed into a direct sum of functions between interval modules. We also give a weaker version of the converse problem, and discuss with examples.

We start by showing that non-trivial maps between interval modules are sent to rectangles (or zero) by $S$.

\begin{lemma}\label{lemma: S_rectangle}
    For intervals $\langle a^\bullet,b^\bullet\rangle$ and $\langle c^\bullet,d^\bullet\rangle$ in $\mathbb{R}$, if $\phi:\kk_{\langle a^\bullet,b^\bullet\rangle}\to\kk_{\langle c^\bullet,d^\bullet\rangle}$ is a non-zero persistence map, then 
    \begin{equation*}  S(\phi) \cong \kk_{ \langle c^\bullet,\infty)\times \gint{a}{d}\times \langle b^\bullet,\infty)}.
    \end{equation*}
\end{lemma}
\begin{proof}
    Let $r$, $t$ and $s$ be real numbers. Pointwise, we have    
    \begin{equation*}
        S(\phi)_{(r,t,s)}=(\im \kk_{\langle c^\bullet,d^\bullet\rangle})_{(r,t)}\cap\phi_t((\ker \kk_{\langle a^\bullet,b^\bullet\rangle})_{(t,s)}).
    \end{equation*}
    From \cref{lemma: image_is_rectangle}, the image $(\im \kk_\gint{c}{d})_{(r,t)}$ is zero if either $r\not\in\langle c^\bullet,\infty)$ or $t\not\in\gint{c}{d}$, and the kernel $(\ker \kk_{\langle a^\bullet,b^\bullet\rangle})_{(t,s)}$, and thus the submodule $\phi_t((\ker \kk_{\langle a^\bullet,b^\bullet\rangle})_{(t,s)})$) is zero if either $t\not\in\langle a^\bullet,b^\bullet\rangle$ or $s\not\in\langle b^\bullet,\infty)$. 
    
    Since $\phi$ is non-zero, we know that $c^\bullet\leq a^\bullet\leq d^\bullet \leq b^\bullet$, so $S(\phi)_{(r,t,s)}$ is only non-zero whenever $r\in\langle c^\bullet,\infty)$, $t\in\gint{a}{d}$ and $s\in \langle b^\bullet,\infty)$, in which case it equals $\im\phi_t\cong \kk$. Structure maps $S(\phi)_{(r,t,s),(r',t',s')}$ are given by restrictions from $(\kk_{\langle c^\bullet,d^\bullet\rangle})_{t,t'}$.
\end{proof}

\begin{remark}
    Notice that $S(\alpha)$ does not detect trivial persistence morphisms, i.e if $\alpha: 0\to \kk_J$ or $\alpha : \kk_I\to 0$, then $S(\alpha)=0$. 
    Also, by definition, if $\ker V=0$, then $S(\alpha)=0$. Hence, $S(\alpha)$ is not affected by elements from $V / \ker(V)$; i.e. generators of $V$ corresponding to intervals from $\B(V)$ with an infinite right endpoint. For instance, considering $\beta: \kk_{[0,\infty)}\to \kk_{[-1,1]}$, it holds that $S(\beta)=0$.
\end{remark}

\begin{example}[$S(\alpha)$ has more information than $\alpha(V)$ and $\ker\alpha$\label{example: alpha_has_less_info_X}]
   Consider the following two ladder modules, for fixed $c\leq b$:
    \begin{align*}
    \alpha&:\kk_{[1,3]}\to \kk_{[0,2]},\\
    \beta_{bc} &:\kk_{[1,b]}\to \kk_{[c,2]}.
    \end{align*}
    Note that their images are isomorphic when $c<1$ and $b>3$, as $\alpha(\kk_{[1,3]})\cong \kk_{[1,2]}\cong \beta_{bc}(\kk_{[1,b]})$. However, computing $S(\alpha)$ and $S(\beta_{bc})$ we see that
\[ S(\alpha)\cong \kk_{[0,\infty)\times [1,2]\times (3,\infty)}\ncong \kk_{[c,\infty)\times [1,2]\times (b,\infty)}\cong S(\beta_{bc}).\]
Hence, $S$ is captures more information about $\alpha$ and $\beta_{bc}$ than the image module alone. Analogously, fixing $b=3$, we get $\ker \alpha \cong  \kk_{(2,3]}\cong \ker \beta_{bc}$ while $S(\alpha)\ncong S(\beta_{bc})$ for $c<1$.
\end{example}

Let $V$ be an interval decomposable module with barcode $\B(V)$. Define $\widehat{\B}(V)\coloneq \B(V)\cup \{(\emptyset, \infty)\}$ where $\{(\emptyset, \infty)\}$ is the multiset containing the empty set with infinite multiplicity. 
Also, recall our convention that $\kk_\emptyset =0$. A ladder module $\alpha:V\to U$ between interval decomposable modules is \emph{bar-to-bar} if it decomposes into a direct sum of the form
\begin{equation*}
    \alpha  \cong \bigoplus_{(I,J)\in \Lambda} (\phi_{IJ}:\kk_{I} \rightarrow \kk_{J})
\end{equation*}
where $\Lambda\subseteq \widehat{\B}(V)\times\widehat{\B}(U)$ and each $\phi_{IJ}$ is an indecomposable map. We write $\Lambda_0\coloneq\Lambda\cap (\B(V)\times \B(U))$ for the pairs corresponding to maps between two non-empty intervals. In particular, $\phi_{IJ}$ is non-zero if and only if $(I,J)\in\Lambda_0$.

\begin{theorem}\label{theorem: bar_to_rect}
    If $\alpha:V\to U$ is bar-to-bar, then $S(\alpha)$ is rectangle decomposable. Specifically, if  $\alpha\cong\bigoplus_{(I,J)\in\Lambda} (\phi_{IJ}:\kk_I\to\kk_J)$ is a bar-to-bar decomposition, then
    \begin{equation*} 
    S(\alpha)\cong \bigoplus_{(\gint{a}{b},\gint{c}{d}) \in \Lambda_0}\kk_{\langle c^\bullet,\infty)\times \gint{a}{d} \times \langle b^\bullet,\infty)},
    \end{equation*}
    where $(I,J)\in\Lambda_0\subseteq \Lambda$ if and only if $\phi_{IJ}\neq 0$.
\end{theorem}

\begin{proof}
    Let $\alpha\cong\bigoplus_{(I,J)\in\Lambda} \phi_{IJ}$ where $\phi_{IJ}:\kk_{I}\to\kk_{J}$. By functoriality (\cref{prop: S_functorial}) and additivity (\cref{prop: S_additive}) of $S$, we get
    \begin{equation*}
         S(\alpha)\cong  S\left(\bigoplus_{(I,J)\in\Lambda} \phi_{IJ}\right)\cong \bigoplus_{(I,J)\in\Lambda} S(\phi_{IJ}),
    \end{equation*}
    and the result follows from \cref{lemma: S_rectangle}, noting that $S(\phi)=0$ when $\phi=0$.
\end{proof}
Again, we stress that if $I=\gint{a}{b}$ with $b=\infty$, then $S(\phi_{IJ})=0$. The theorem is true as stated, remembering the convention $\kk_\emptyset=0$.

\begin{example}[Computation of $S(\alpha)$ for $\alpha$ indecomposable]\label{example:Scomputation}
Consider the persistence modules $V=\kk_{[2,4]}\oplus \kk_{[1,5]}$ and $U=\kk_{[0,3]}$ with the ladder module
\[ \alpha:\kk_{[2,4]}\oplus \kk_{[1,5]} \xrightarrow{[1\ 1]}  \kk_{[0,3]},\]
mapping the two bars of the domain onto the same bar in the target. We compute $S(\alpha)$, first using \cref{corollary: imkerdecomp} to compute (see \cref{fig: S_computation_example})
\[ \ker V\cong \kk_{[2,4]\times (4,\infty)} \oplus \kk_{[1,5]\times(5,\infty)}  \]

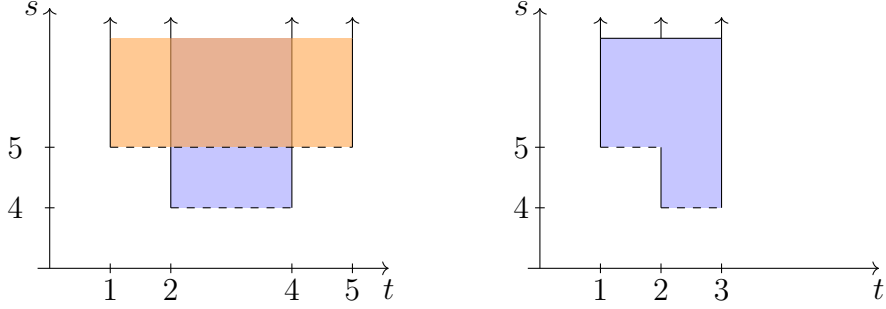
\begin{figure}[ht]
\centering

\begin{subfigure}{0.48\textwidth}
\centering
\begin{tikzpicture}[scale=0.8]

\draw[->] (-0.2,0) -- (5.6,0) node[below] {$t$};
\draw[->] (0,0) -- (0,4.3) node[left] {$s$};

\fill[blue!50,opacity=.45]
    (2,1) rectangle (4,3.8);

\draw[dashed] (2,1) -- (4,1);
\draw (2,1) -- (2,3.8);
\draw (4,1) -- (4,3.8);

\fill[orange!70,opacity=.60]
    (1,2) rectangle (5,3.8);

\draw[dashed] (1,2) -- (5,2);
\draw (1,2) -- (1,3.8);
\draw (5,2) -- (5,3.8);

\draw[->] (2,3.8) -- (2,4.15);
\draw[->] (4,3.8) -- (4,4.15);
\draw[->] (1,3.8) -- (1,4.15);
\draw[->] (5,3.8) -- (5,4.15);

\foreach \x in {1,2,4,5}
    \draw (\x,0.08) -- (\x,-0.08);

\foreach \y/\lab in {1/4,2/5}
    \draw (0.08,\y) -- (-0.08,\y)
          node[left=4pt] {$\lab$};

\node[below] at (1,0) {$1$};
\node[below] at (2,0) {$2$};
\node[below] at (4,0) {$4$};
\node[below] at (5,0) {$5$};

\end{tikzpicture}
\end{subfigure}
\hfill
\begin{subfigure}{0.48\textwidth}
\centering
\begin{tikzpicture}[scale=0.8]

\draw[->] (-0.2,0) -- (5.6,0) node[below] {$t$};
\draw[->] (0,0) -- (0,4.3) node[left] {$s$};

\fill[blue!50,opacity=.45]
    (1,2) rectangle (3,3.8);

\fill[blue!50,opacity=.45]
    (2,1) rectangle (3,2);

\draw[dashed] (1,2) -- (2,2);
\draw[dashed] (2,1) -- (3,1);

\draw (1,2) -- (1,3.8);
\draw (1,3.8) -- (3,3.8);
\draw (3,3.8) -- (3,1);
\draw (2,1) -- (2,2);

\draw[->] (1,3.8) -- (1,4.15);
\draw[->] (2,3.8) -- (2,4.15);
\draw[->] (3,3.8) -- (3,4.15);

\foreach \x in {1,2,3}
    \draw (\x,0.08) -- (\x,-0.08);

\foreach \y/\lab in {1/4,2/5}
    \draw (0.08,\y) -- (-0.08,\y);

\node[below] at (1,0) {$1$};
\node[below] at (2,0) {$2$};
\node[below] at (3,0) {$3$};

\node[left] at (0,1) {$4$};
\node[left] at (0,2) {$5$};

\end{tikzpicture}
\end{subfigure}

\caption{ Figures of example \ref{example:Scomputation} \underline{Left}: Rectangles in the decomposition of $\ker V$. \underline{Right}: The region $\Lambda$, where $\kk_\Lambda=\alpha(\ker V)$.}
\label{fig: S_computation_example}
\end{figure}

Since the two bars of $V$ are merged in the image, and that for every $t>d$ and every $s$ we have $\alpha_t(\ker V_{t,s})=0$, we get that the indicator module $\kk_\Lambda$ where $\Lambda=([1,2)\times(5,\infty))\cup([2,3]\times(4,\infty))$ is shown in \cref{fig: S_computation_example}.

To finish the computation of $S(\alpha)$, we note that $\im U\cong \kk_{[0,\infty]\times [0,3)}$ by \cref{corollary: imkerdecomp}, and so intersection $(\im U)_{r,t}$ and $\alpha_t(\ker V)_{t,s}$ we get $S(\alpha)\cong\kk_{[0,\infty)\times \Lambda}$ which is not decomposable.
\end{example}

\begin{example}[Infinite bars]\label{example: problems_with_infinities}
    Consider the persistence modules $V=\kk_{[1,\infty)}\oplus \kk_{[2,4]}$ and $U=\kk_{[0,3]}$ with the indecomposable ladder module
    \[ \alpha:\kk_{[1,\infty)}\oplus \kk_{[2,4]} \xrightarrow{[1\ 1]}  \kk_{[0,3]}.\]

    Since the kernel sends infinite bars to zero, we get $\ker V \cong \kk_{[2,4]\times (4,\infty)}$, and combining it with $\im U\cong \kk_{[0,\infty)\times[0,3]}$, we compute
    \[ S(\alpha)\cong \kk_{[0,\infty]\times [2,3]\times (4,\infty)}.\]

\end{example}
\cref{example: problems_with_infinities} shows a ladder module $\alpha$ that is not bar-to-bar, but where $S(\alpha)$ still is a rectangle, giving a counterexample to the converse of \cref{theorem: bar_to_rect}. In the rest of this section, we will work towards a weaker version of this converse statement. We start by showing that rectangle summands in $S(\alpha)$ must have the shape found in \cref{lemma: S_rectangle}. 

\begin{proposition}\label{prop: shape-rectangles}
If $\alpha:V\to U$ is a ladder module where $S(\alpha)\cong \kk_R \oplus A$ for a rectangle $R$, then $R={\langle c^\bullet,\infty)\times \langle a^\bullet,d^\bullet \rangle\times\langle b^\bullet,\infty)}$ where $c^\bullet\leq a^\bullet\leq d^\bullet\leq b^\bullet$.
\end{proposition}
\begin{proof}
    Suppose that $R=\gint{c}{x}\times\gint{a}{d}\times\gint{b}{y}$. We have $x^\bullet=y^\bullet=\infty^\bullet$ since structure maps of $S(\alpha)$ are the identity in the first and third parameters. We also have $a^\bullet\leq d^\bullet$ since $\gint{a}{d}$ is non-empty. 

    Assume by contradiction that $b^\bullet<d^\bullet$, so there is a point $(r,t,s)\in R$ with $r\in \langle c^\bullet,\infty) $ and $t>s$. Since $\ker V_{t,s}=0$, this gives a contradiction: 
    \[0\neq \kk\oplus A_{(r,t,s)}=(\kk_R \oplus A)_{(r,t,s)}\cong S(\alpha)_{(r,t,s)}=\im U_{r,t}\cap \alpha_t(\ker V_{t,s})=0.\]

    Similarly, if we assume $a^\bullet<c^\bullet$, and let $t\in \langle a^\bullet,c^\bullet \rangle$, $r\in \langle c^\bullet,\infty)$ and $s\in \langle b^\bullet,\infty)$. Then $(r,t,s)\in R$ while $(t,t,s)\not\in R$, and  we have
    \[
    S(\alpha)_{(t,t,s)}=\alpha_t(\ker V_{t,s})=S(\alpha)_{(r,t,s)}\,.
    \]
    However, restricting to $\kk_R$, this implies that
    \begin{equation*}
        0=(\kk_R)_{(t,t,s)}=(\kk_R)_{(r,t,s)}=\kk
    \end{equation*}
    leading to a contradiction.
\end{proof}

Since the $S$ functor is blind to infinite bars in the domain, we choose to ignore them by restricting to persistence modules with only bars of finite length; see \Cref{def:isolated-critical-finite-length}. In order to tame the case in which our modules have an infinite number of bars, we also consider modules with an isolated set of critical points; see \Cref{def:isolated-critical-finite-length}.

\begin{restatable}{theorem}{converseresult} \label{thm: Rect_to_bar}
    Let $\alpha: V\to U$ be a ladder module with isolated critical points, where $V$ has bars of finite length. If $S(\alpha)\cong \kk_R$ for some rectangle $R$, then $\alpha$ is bar-to-bar with a single non-trivial summand. Furthermore, if $R={\langle c^\bullet,\infty)\times \langle a^\bullet,d^\bullet \rangle\times\langle b^\bullet,\infty)}$, then the summand is $\kk_{\gint{a}{b}}\xrightarrow{}\kk_\gint{c}{d}$.
\end{restatable}
\begin{proof}
    The proof is given in \cref{appendixB}.
\end{proof}

As a summary, we have defined $S$, an assignment between ladder modules and 3-parameter persistence modules, and showed that it is \textit{ functorial, stable with respect to the interleaving distance, additive} and that it sends bar-to-bar ladder modules into rectangle decomposable modules. Furthermore, \cref{thm: Rect_to_bar} suggests that rectangle summands in $S$ hold information about the individual bar-to-bar components of $\alpha$.

\section{Conclusion}

In this work, we presented the additive functors $\ker(-)$ and $\im(-)$, which turn out to preserve the interleaving distance and encode the life and death of features in one-parameter persistence modules. These two functors send an interval-decomposable one-parameter persistence module (with finite bars for the $\ker(-)$ functor) to a rectangle decomposable two-parameter persistence module. They are the building blocks of our new invariant for ladder modules, given by the functor $S$.

Throughout the paper, we proved that the functor $S$ is stable with respect to the interleaving distance and is additive. We also proved that $S$ sends bar-to-bar morphisms to rectangle-decomposable modules. Furthermore, we showed that if $S(\alpha)$ is a rectangle, then it must have a specific shape, and that $\alpha$ must be the direct sum of a morphism between two intervals and trivial morphisms.

Our new invariant can be seen as a different version of additive matching, in which we trade some of the combinatorial simplicity for stability and functoriality. In this way, our work connects the theory of ladder modules with the theory of rectangle-decomposable modules in multipersistence.

As future work, we would like to develop an algorithm to compute the functor $S$ for a ladder module $\alpha$ under mild conditions, and to connect it with the theory of additive partial matchings~\cite{GONZALEZDIAZ2026102598,Guidolin-matchings}. We also believe that some of the assumptions we impose in the proof of the converse result \ref{thm: Rect_to_bar} are too strong. We conjecture that the same result holds for rectangle-decomposable modules, and we are currently working on this direction.

In addition, in this work we have introduced a general notion of interleaving between multi-parameter persistence modules. In particular, we would like to relate \Cref{thm: shift_respecting_stable} and \Cref{def: respecting_shifts} with existing results on general interleavings, such as Theorem~4.3 from~\cite{categorieswithaflow}.

\bibliographystyle{plainurl}
\bibliography{bibliography}

\appendix
\crefalias{section}{appendix}

\section{Stability of shift respecting functors:}\label{sec: AppendixA}

\begin{theorem}
    Let $\Sigma^1$ and $\Sigma^2$ be shifts of $P_1$ and $P_2$ respectively, and let $V$ and $U$ be persistence modules over $P_1$. If $F:\veck^{P_1}\to\veck^{P_2}$ respects $\Sigma^1$ and $\Sigma^2$, then 
    \begin{equation*}
        d^{\Sigma^1}_I(V,U)\geq d^{\Sigma^2}_I(F(V),F(U)).
    \end{equation*}
\end{theorem}

\begin{proof}
    Let $(\phi,\psi)$ be an $\varepsilon$-interleaving between $V$ and $U$, so we have morphisms $\phi:V\to U[\varepsilon]$ and $\psi:U\to V[\varepsilon]$ where $\psi[\varepsilon]_p\circ \phi_p = V_{p,\Sigma^1_{2\varepsilon}(p)}$ and $\phi[\varepsilon]_p\circ \psi_p = U_{p,\Sigma^1_{2\varepsilon}(p)}$ for all $p\in P_1$. 
    
    Since $F$ respects shifts, there is a collection of natural maps $N^\varepsilon:F((-)[\varepsilon])\to F(-)[\varepsilon]$ for $\varepsilon\geq 0$. Consider $\Phi\coloneq N^{\varepsilon}_U\circ F(\phi)$ and $\Psi\coloneq N^{\varepsilon}_V\circ F(\psi)$. By functoriality of the shifting, we get
    \begin{equation*}
        \Phi[\varepsilon] = N^{\varepsilon}_{U}[\varepsilon]\circ F(\phi)[\varepsilon],
    \end{equation*}
    and similarly, $\Psi[\varepsilon] = N^{\varepsilon}_V[\varepsilon] \circ F(\psi)[\varepsilon]$.

    Consider the composition  
    \[ F(V)\xrightarrow[]{\Phi}F(U)[\varepsilon]\xrightarrow[]{\Psi[\varepsilon]}F(V)[2\varepsilon].\]
Using the various compatibility conditions of \cref{def: respecting_shifts} together with the fact that $(\phi,\psi)$ is an interleaving and $F$ is a functor, we get
\begin{align*}
\Psi[\varepsilon]\circ \Phi&=N_V^{\varepsilon}[\varepsilon]\circ F(\psi)[\varepsilon]\circ N_U^{\varepsilon}\circ F(\phi)
\\&=N^{\varepsilon}_V[\varepsilon ]\circ N_{V[\varepsilon]}^{\varepsilon}\circ F(\psi[\varepsilon])\circ F(\phi)\\&=N^{2\varepsilon}_V\circ  F(\psi[\varepsilon])\circ F(\phi)\circ(N^0_V)^{-1}\\&=N^{2\varepsilon}_V\circ F(\tau_{0,2\varepsilon}^{V,\Sigma^1})\circ(N^0_V)^{-1}=\tau_{0,2\varepsilon}^{F(V),\Sigma^2}
\end{align*}

By symmetry, we get that $(\Phi,\Psi)$ is an $\varepsilon$-interleaving between $F(V)$ and $F(U)$ (with respect to $\Sigma^2$).
\end{proof}

\begin{remark}
The commutativity of the following diagram summarizes the proof:
{\small
\[
\hspace{-1cm}\begin{tikzcd}[sep=scriptsize]
	& {F(V)} && {F(V)[\varepsilon]} && {F(V)[2\varepsilon]} & \\
	& {} & {F(U)} & {} & {F(U)[\varepsilon]} && {F(U)[2\varepsilon]} \\
	{} & {} \\
	& {F(V)} && {F(V[\varepsilon])} && {F(V[2\varepsilon])} \\
	&& {F(U)} && {F(U[\varepsilon])} && {F(V[2\varepsilon])}
	\arrow[from=1-2, to=1-4]
	\arrow["\Phi"{description, pos=0.3}, from=1-2, to=2-5]
	\arrow[from=1-4, to=1-6]
	\arrow["{\Phi[\varepsilon]}"{description, pos=0.3}, from=1-4, to=2-7]
	\arrow[from=2-3, to=1-4]
	\arrow[from=2-3, to=2-5]
	\arrow[from=2-5, to=1-6]
	\arrow[from=2-5, to=2-7]
	\arrow["{id_{F(V)}}"', from=4-2, to=1-2]
	\arrow[from=4-2, to=4-4]
	\arrow["{F(\phi)}"{description, pos=0.8}, from=4-2, to=5-5]
	\arrow["{N^{\varepsilon}_V}"', from=4-4, to=1-4]
	\arrow[from=4-4, to=4-6]
	\arrow["{F(\phi[\varepsilon])}"{description, pos=0.7}, from=4-4, to=5-7]
	\arrow["{N^{2\varepsilon}_V}"', from=4-6, to=1-6]
	\arrow["{id_{F(U)}}"', from=5-3, to=2-3]
	\arrow[from=5-3, to=4-4]
	\arrow[from=5-3, to=5-5]
	\arrow["{N^{\varepsilon}_U}"', from=5-5, to=2-5]
	\arrow[from=5-5, to=4-6]
	\arrow[from=5-5, to=5-7]
	\arrow["{N^{2\varepsilon}_U}"', from=5-7, to=2-7]
\end{tikzcd}\]}
The intuition behind the definition of a shift respecting functor is making every rectangle in this diagram commute. 
\end{remark}
\begin{remark}\label{equality case}
    Let $P_1$ and $P_2$ be posets with shifts $\Sigma^1$ and $\Sigma^2$ respectively. If $F:\veck^{P_1}\to \veck^{P_2}$ is a functor such that for all $\varepsilon\geq 0$ and $V\in \veck^{P_1}$,
    \begin{itemize}
        \item[(1)]$F(V[\varepsilon])= F(V)[\varepsilon]$,
      
        \item[(2)] $F(\tau_{\varepsilon,\varepsilon'}^{V,\Sigma^1})=\tau_{\varepsilon,\varepsilon'}^{F(V),\Sigma^2}$ for $\varepsilon\leq \varepsilon'$.
    \end{itemize}
      Then, it is a direct consequence of \Cref{def: respecting_shifts} that $F$ respects the shifts.
\end{remark}

\begin{lemma}\label{lemma: shifting_X}
    Let $M$ be a $\mathbb{R}$-persistence module and $\alpha:V\to U$ be a ladder module over $\mathbb{R}$. If $\varepsilon\geq 0$, then:
    \begin{itemize}
        \item $\im(M[\varepsilon])=(\im M)[\varepsilon]$
        \item $\ker(M[\varepsilon])=(\ker M)[\varepsilon]$
        \item $S(\alpha[\varepsilon])=S(\alpha)[\varepsilon]$.
    \end{itemize}
\end{lemma}
\begin{proof}
   Note that for a persistence module $M$ and $\varepsilon>0$ we have
\begin{equation*}
    (\im M[\varepsilon])_{r,t}=
    \begin{cases}
        \im M[\varepsilon]_{(r,t)} = \im M_{r+\varepsilon,t+\varepsilon} &\textrm{if }r\leq t,\\
        M[\varepsilon]_t = M_{t+\varepsilon} & \textrm{if }r> t,
    \end{cases}
\end{equation*}
which is exactly $(\im M)_{r+\varepsilon,t+\varepsilon}=(\im M)[\varepsilon]_{(r,t)}$. In particular, $\im M[\varepsilon] = (\im M)[\varepsilon]$, and similarly we have $\ker M[\varepsilon] = (\ker M)[\varepsilon]$.

To proof the equality for the funtor $S$, recall that $\alpha[\varepsilon]:V[\varepsilon]\to U[\varepsilon]$, where $V[\varepsilon]_t=V_{t+\varepsilon}$ and $V[\varepsilon]_{t,t'}=V_{t+\varepsilon,t'+\varepsilon}$, and similar for $U$, and where $\alpha[\varepsilon]_t=\alpha_{t+\varepsilon}$. In particular, we get
    \begin{align*}
        S(\alpha[\varepsilon])_{(a,t,b)}&= (\im U)_{(a+\varepsilon,t+\varepsilon)}\cap\alpha_{t+\varepsilon}((\ker V)_{(t+\varepsilon,b+\varepsilon)}) \\
        &= S(\alpha)_{(a+\varepsilon,t+\varepsilon,b+\varepsilon)} = S(\alpha)[\varepsilon]_{(a,t,b)}.
    \end{align*}
    concluding the result.
\end{proof}
\begin{lemma}\label{lemma: S_im_ker_respect}
     The functors $S$, $\im$ and $\ker$ are shift respecting for the shifts defined in \ref{sec: interleaving distance}.
\end{lemma}
\begin{proof}
    Applying the argument from Remark $\ref{equality case}$ and using Lemma \ref{lemma: shifting_X}, we just need to check that the functors $S$, $\im$ and $\ker$ preserve the canonical shifting map in the sense of Remark \ref{equality case}. By definition of $S$: 
   \begin{align*}
    S(\tau_{\varepsilon,\varepsilon'}^{\alpha,\Sigma^1})_{(r,t,s)}&=(U_{t+\varepsilon,t+\varepsilon'})_{\mid S(\alpha)_{(r+\varepsilon,t+\varepsilon,s+\varepsilon)}}\\&=S(\alpha)_{((r+\varepsilon,t+\varepsilon,s+\varepsilon),(r+\varepsilon',t+\varepsilon',s+\varepsilon'))}=(\tau_{\varepsilon,\varepsilon'}^{S(\alpha),\Sigma_2})_{(r,t,s)}
    \end{align*}
    For the functors $\im$ and $\ker$ a similar arguments also works.
\end{proof}

\begin{lemma}\label{lemma: imkerstab}
For persistence modules $V$ and $U$ over $\mathbb{R}$, we have
   \[ \dint(V,U) \geq \dintn{2}(\im V ,\im U )
   \qquad\textrm{and}\qquad
    \dint(V,U) \geq \dintn{2}(\ker V ,\ker U ).
   \]
\end{lemma}
\begin{proof}
Applying Lemma \ref{lemma: S_im_ker_respect} and Theorem \ref{thm: shift_respecting_stable} we get the result.
  
\end{proof}

Even more, we can prove the converse inequality, assuming no infinite bars in the case of the kernel:

\begin{lemma}\label{lemma: imkerstab2}
    For persistence modules $V$ and $U$ over $\mathbb{R}$, we have
    \[\dint(V,U)\leq\dintn{2}(\im V ,\im U ).\]
    Furthermore, if $V$ and $U$ have bars of finite length, then
    \[ \dint(V,U)\leq\dintn{2}(\ker V ,\ker U ).\]
\end{lemma}

\begin{proof}
\textit{Image:} Note that $(\im V)_{t,t}=V_t$, so if $\phi:\im V\to (\im U)[\varepsilon]$ and $\psi: \im U\to (\im V)[\varepsilon]$ form an $\varepsilon$-interleaving, then the diagonal maps $(\phi_{(-,-)},\psi_{(-,-)})$ gives an $\varepsilon$-interleaving between $V$ and $U$.

\textit{Kernel:} Since $V$ has bars of finite length in its barcode, for each $t\in \mathbb{R}$, every element $v\in V_t$ lies in $\ker V_{t,s}$ for some $s\geq t$. So taking the colimit, we get an isomorphism:
\[ \varinjlim_{s\in \mathbb{R}}\left(\ker V_{-,s}\right)\cong V.\]
We can also see the the functor:
\[ 
\varinjlim_{s\in \mathbb{R}}(-):\veck^{\mathbb{R}^2}\to \veck^{\mathbb{R}}
\]
sending $M\mapsto \varinjlim_{s\in \mathbb{R}}(M_{-,s})$ is shift respecting for the shifts defined in \Cref{sec: interleaving distance}.
Hence, applying \Cref{thm: shift_respecting_stable} we get:
\[
d_I^{\mathbb{R}^2}(\ker V,\ker U)\geq d_I^{\mathbb{R}}\left(\varinjlim_{s\in \mathbb{R}}(\ker V_{-,s}), \varinjlim_{s\in \mathbb{R}}(\ker U_{-,s}) \right)=\dint(V,U).
\]
\end{proof}

\isoimker*
\begin{proof}
    The result is a direct consequence of Lemma \ref{lemma: imkerstab}  and Lemma \ref{lemma: imkerstab2}
\end{proof}

We end up Appendix A proving the last result from Section \ref{sec: Stability}:

\begin{proposition}
    Given a pair of ladder modules $\alpha:V\to U , \alpha':V'\to U'$, if $V$ and $V'$ have bars of finite length, then
    \[  \dint (\alpha(V),\alpha'(V')) \leq d_I^{\mathbb{R}^3}(S(\alpha),S(\alpha')).\]
\end{proposition}
\begin{proof}
    If $V$ has bars with finite length we have:

\[\underset{(r,s)\in \mathbb{R}^2}{\varinjlim}S(\alpha)_{(r,-,s)}=\underset{(r,s)\in \mathbb{R}^2}{\varinjlim} (\im U_{r,-})\cap \alpha(\ker V_{-,s}) \cong \alpha(V)\]
 Since that colimit functor is shift preserving, we can apply Theorem \ref{thm: shift_respecting_stable} to get:
 \begin{align*}
 d_I^{\mathbb{R}^3}(S(\alpha),S(\alpha')) & \geq \dint\left(\underset{(r,s)\in \mathbb{R}^2}{\varinjlim}S(\alpha)_{(r,-,s)},\underset{(r,s)\in \mathbb{R}^2}{\varinjlim}S(\alpha')_{(r,-,s)}\right) \\
    & =\dint (\alpha(V),\alpha'(V'))\,.
 \hfill \qedhere
 \end{align*}
\end{proof}

\section{Proof of \texorpdfstring{\cref{thm: Rect_to_bar}}{Theorem}}\label{appendixB}

In this section, we prove the following:
\converseresult*

Since $\kk_R$ is a summand of $S(\alpha)$, there exist persistence maps
\begin{equation}\label{eq: Psi_definition}
    \kk_R\xhookrightarrow{\Psi} S(\alpha)\twoheadrightarrow\kk_R
\end{equation}
whose composition is the identity. 
For $(r,t,s)\in R$, we define $u_t\coloneq \Psi_{(r,t,s)}(1_\kk)$, where $1_\kk$ is the identity in $(\kk_R)_{(r,t,s)}=\kk$. 
Note that $u_t$ is non-zero, and that
\begin{equation}\label{eq: ut_in_im_and_ker}
    u_t\in S(\alpha)_{(r,t,s)}=\im(U_{r,t}) \cap \alpha_t(\ker(V_{t,s})).
\end{equation}

\begin{remark}
    Indeed $u_t$ is well-defined: If $x\coloneq(r,t,s) \in R$ and $x'=(r',t,s')\in R$, then $m\coloneq(\min\{r,r'\},t,\min\{s,s'\})$ is in $R$ with $m\leq x,x'$. Using naturality of $\Psi$, we have $\Psi_x \circ (\kk_R)_{m,x}= S(\alpha)_{m,x} \circ \Psi_m$.
    The structure maps $S(\alpha)_{m,x}=U_{t,t}$ and $(\kk_R)_{m,x}$ are identities, so we get $\Psi_m(1_\kk)=\Psi_x(1_\kk)=u_t$. Similarly we have $\Psi_{x'}(1_\kk)=\Psi_m(1_\kk)$.
\end{remark}

Using \cref{prop: shape-rectangles}, we can write $R=\langle c^\bullet, \infty)\times \gint{a}{d}\times\langle b^\bullet,\infty)$.
Since $\alpha$ has isolated critical points, we can find $\wa\in\gint{a}{d}$, $\wb\in\langle b^\bullet ,\infty)$ and $\wc \in \gint{c}{d}$ making $\alpha$ constant along $\langle a^\bullet, \wa]$, $\langle b^\bullet,\wb]$ and $\langle c^\bullet, \wc]$ respectively. Without loss of generality we assume that $\wc \leq \wa$.

By \cref{eq: ut_in_im_and_ker}, we have $u_{\wa}\in \im(U_{\wc,\wa})$, so let $v \in U_{\wc}$ be such that $U_{\wc,\wa}(v)=u_\wa\neq 0$, implying $v\neq 0$.
Define the map $\varepsilon\colon \kk_{\gint{c}{d}}\rightarrow U$ for $t \in \gint{c}{d}$ by
\begin{equation} \label{eq: varepsilon_def}
    \varepsilon_t(1_\kk)=\begin{cases}
    v & \mbox{ for } t \in \langle c^\bullet,\wc] \\
    U_{\wc,t}(v) & \mbox{ for } t \in [\wc,d^\bullet\rangle.
\end{cases}
\end{equation}
This is a well-defined persistence map, since $U_t=U_\wc$ for $t \in \langle c^\bullet,\wc]$ and $\varepsilon_t(1_\kk)=U_{t,t'}(\varepsilon_{t'}(1_\kk))$ for all $t\leq t'$ in $\gint{c}{d}$.  
In particular, $\varepsilon_t(1_\kk)=u_t$ whenever $t\in\gint{a}{d}$, because $U_{\wc,t}(v)=u_\wa=u_t$ for $t\in \langle a^\bullet, \wa]$ and $U_{\wc,t}(v)=U_{\wa,t}(u_\wa)=U_{\wa,t}(\Psi_{(r,\wa,s)}(1_\kk))=\Psi_{(r,t,s)}(1_\kk)=u_t$ for $t\in [\wa,d^\bullet\rangle$.

Also by \cref{eq: ut_in_im_and_ker}, we have $u_\wa\in \alpha_\wa(\ker V_{\wa,\wb})$, so there exists $w\in\ker(V_{\wa, \wb})$ such that $\alpha_{\wa}(w)=u_{\wa}$. Fixing such a $w$, define the map $\delta\colon \kk_{\gint{a}{b}}\rightarrow V$ for $t\in\gint{a}{b}$ by
\begin{equation} \label{eq: delta_def}
    \delta_t(1_\kk)=
\begin{cases}
    w & \mbox{ for } t \in \langle a^\bullet, \wa] \\
    V_{\wa, t}(w) & \mbox{ for } t \in [\wa, b^\bullet\rangle\,.
\end{cases}
\end{equation}
Altogether, we obtain the following natural transformation
\begin{equation}\label{diag:commutativity-alpha-one}
 \begin{tikzcd}[/tikz/column 2/.style={column sep=0cm}]
    \kk_{\gint{a}{b}} \ar[r, "1"] \ar[d, "\delta"] &
    \kk_{\gint{c}{d}} \ar[d, "\varepsilon"] \\
    V \ar[r, "\alpha"] & U 
    & .
\end{tikzcd}   
\end{equation}

\begin{lemma}\label{lemma: properties-delta-epsilon}
    Let $\alpha: V\to U$ be a ladder module with isolated critical points such that \(S(\alpha)\cong \kk_{\langle c^\bullet,\infty)\times \gint{a}{d}\times\langle b^\bullet,\infty)}\oplus A\).
    If the maps $\varepsilon\colon \kk_{\gint{c}{d}}\to U$ and $\delta \colon \kk_{\gint{a}{b}} \to V$ are as in \cref{eq: varepsilon_def,eq: delta_def}, then 
\begin{enumerate}[label=(\roman*)]
    \item \label{delta_epsilon:naturality}$\alpha_t(\delta_t(1_\kk))=\varepsilon_t(1_\kk)$ for all $t \in \gint{a}{d}$,
    \item\label{epsilon:non-zero} $\varepsilon_t\neq 0$ for all $t\in \langle c^\bullet, d^\bullet\rangle$,
    \item\label{delta:nan-zero} $\delta_t\neq 0$ for all $t\in \langle a^\bullet, b^\bullet\rangle$,
    \item $\im U_{r,t}\cap\im \varepsilon_t = 0$ for all $r \in \langle -\infty, c^\bullet\rangle$ and all $t \in \langle c^\bullet, d^\bullet\rangle$.
\end{enumerate}
Furthermore, if $\alpha_r=0$ for all $r\in(-\infty,a^\bullet\rangle$, then
\begin{enumerate}[label=(\roman*)]
    \setcounter{enumi}{4}
    \item $\im V_{r,t}\cap\im \delta_t = 0$ for all $r \in \langle -\infty, a^\bullet\rangle$ and all $t \in \langle a^\bullet, b^\bullet\rangle$.
\end{enumerate}
\end{lemma}
\begin{proof}
    Write $R={\langle c^\bullet,\infty)\times \langle a^\bullet,d^\bullet \rangle\times\langle b^\bullet,\infty)}$. We show each property:
    \begin{enumerate}[label=(\roman*)]
        \item This follows by the naturality of diagram~\eqref{diag:commutativity-alpha-one}.

        \item After \cref{eq: varepsilon_def}, we noted that
        $\varepsilon_t(1_\kk)=u_t$ for all $t\in\gint{a}{d}$.
        This implies that $\varepsilon_t(1_\kk)\neq 0$ for all $t\in\gint{c}{d}$, since, writing $\wt=\max(\wa,t)$, we have $U_{t,\wt}(\varepsilon_t(1_\kk))=\varepsilon_\wt(1_\kk)=u_\wt\neq 0$.
        
        \item By contradiction, assume there exists $t \in \gint{a}{b}$ such that $\delta_t=0$. 
        Using properties \ref{delta_epsilon:naturality} and \ref{epsilon:non-zero}, it follows that $t \in \gint{d}{b}$ since otherwise, we have $t\in\gint{a}{d}$ which implies $\delta_t(1_\kk)=\alpha_t(\varepsilon_t(1_\kk))=u_t\neq 0$.
        With $\wa$ as above, 
        we have $V_{\wa,t}(\delta_\wa(1_\kk))=\delta_t(1_\kk)=0$. In particular, $u_\wa\in\alpha_\wa(\ker(V_{\wa,t}))=S(\alpha)_{\wa,\wa,t}$ where $(\wa,\wa,t)\not\in R$. However, $u_{\wa}$ is sent by structure maps to $u_\wa\in\im \Psi_{(\wa,\wa,s)}$ for $s\in\langle b^\bullet,\infty)$, contradicting the hypothesis that $\im\Psi$ is a direct summand of $S(\alpha)$.

        \item By contradiction, assume there exist $r \in \langle -\infty, c^\bullet\rangle$ and $t \in \langle c^\bullet, d^\bullet\rangle$ such that $\im U_{r,t}\cap\im \varepsilon_t \neq 0$. 
        So, there exists $x \in U_r$ satisfying $U_{r,t}(x)=\varepsilon_t(1_\kk)$. 
        Let $s\in\langle b^\bullet,\infty)$, and note that $V_{t,s}(\delta_t(1_\kk))=\delta_s(1_\kk)=0$ by naturality of $\delta$, so $\delta_{t}(1_\kk) \in \ker(V_{t, s})$. 
        By property \ref{delta_epsilon:naturality}, we have $\varepsilon_{t}(1_\kk) =\alpha_{t}(\delta_{t}(1_\kk)) \in \alpha_{t}(\ker(V_{t, s}))$, and by hypothesis, $\varepsilon_{t}(1_\kk)$ is in $\im (U_{r,t})$, so $\varepsilon_{t}(1_\kk)\in S(\alpha)_{(r,t,s)}$. 
        Note that $(r,t,s)\not\in R$, and if $q\in \gint{a}{d}$ is such that $t\leq q$, then $(t,q,s)\in R$. However, the structure map sends $\varepsilon_t(1_\kk)\not\in \im \Psi_{(r,t,s)}$ to $\varepsilon_q(1_\kk)\in \im \Psi_{(t,q,s)}$ (non-zero by property \ref{epsilon:non-zero}), contradicting the hypothesis that $\im\Psi$ is a direct summand of $S(\alpha)$.

        \item 
        
        By contradiction, assume $r \in \langle -\infty, a^\bullet\rangle$ and $t \in \langle a^\bullet, b^\bullet\rangle$ such that $\im V_{r,t}\cap\im \delta_t \neq 0$, and let $y \in V_r$ satisfy $V_{r,t}(y)=\delta_t(1_\kk)$. 
        For any $q\in \gint{a}{d}$ where $q\leq t$, we have
        \begin{equation*}
            V_{q,t}(\delta_q(1_\kk) - V_{r,q}(y))=\delta_t(1_\kk)-\delta_t(1_\kk)=0,
        \end{equation*}
        and, since $\alpha_r=0$, it follows $\alpha_q(V_{r,q}(y))=U_{r,q}(\alpha_r(y))=0$. Altogether, 
        \[
        0\neq \varepsilon_q(1_\kk) = \alpha_q(\delta_q(1_\kk) - V_{r,q}(y)) \in \alpha_q(\ker(V_{q,t}))=S(\alpha)_{(q,q,t).}
        \]
        Since $(q,q,t)\notin R$, we have $\varepsilon_q(1_\kk)\not\in\im\Psi_{(q,q,t)}$, but it is sent to $\varepsilon_q(1_\kk)\in \im\Psi_{(q,q,s)}$ for any $s\in\langle b^\bullet,\infty)$, contradicting $\im\Psi$ being a direct summand of $S(\alpha).$ \qedhere
    \end{enumerate}
\end{proof}

We continue to show that the interval modules $\im\delta\cong \kk_{\langle a^\bullet,b^\bullet\rangle}$ and $\im\varepsilon\cong \kk_{\langle c^\bullet,d^\bullet\rangle}$ are direct summands of $V$ and $U$, respectively.
\begin{lemma}\label{lemma: phi_splits}
    Let $V$ have isolated critical points, and assume it decomposes into intervals by $\beta: V\xrightarrow\textstyle\bigoplus_{L \in \B(V)}\kk_{L}$. Let $I=\gint{a}{b}$ be an interval. If $\phi \colon \kk_{I} \hookrightarrow V$ is a persistence map with $\phi_t\neq 0$ and $\im(V_{r,t})\cap \im(\phi_t)=0$ for all $r \in (-\infty, a^\bullet\rangle$ and $t \in I$.
    Then, $I\in \B(V)$ and the composition
    \[
\kk_{I} \xhookrightarrow{\phi}
V \xrightarrow[\cong]{\beta} \underset{L\in \B(V) }{\bigoplus}\kk_{L} \twoheadrightarrow \kk_I
\]
is an isomorphism. In particular, $V \cong \kk_{I} \oplus V'$ for some submodule $V'\subset V$.
\end{lemma}
\begin{proof}
For intervals $J \in \B(V)$, we denote the composition
\[
\pi_J \colon \kk_{I} \xhookrightarrow{\phi}
V \xrightarrow[\cong]{\beta} \underset{L\in \B(V) }{\bigoplus}\kk_{L} \twoheadrightarrow \kk_J.
\]
Define the non-empty multiset $\Gamma=\{ J \mid \pi_J\neq 0\}\subseteq \B(V)$. We will show that $\Gamma$ contains a copy of the interval $I=\gint{a}{b}$.

By contradiction, suppose that no copy of $I$ exists in $\Gamma$.
If $J=\gint{c}{d}$ is in $\Gamma$, then the persistence map $\pi_J\colon \kk_{\gint{a}{b}}\rightarrow \kk_{\gint{c}{d}}$ is non-zero, and so $c^\bullet \leq a^\bullet\leq d^\bullet \leq b^\bullet$. 
In particular, we can write $\Gamma = \Gamma_1 \sqcup \Gamma_2$ by setting
\[
\Gamma_1 = \{ \gint{c}{d} \in \Gamma \mid c^\bullet = a^\bullet\}
\mbox{ and }
\Gamma_2 = \{ \gint{c}{d} \in \Gamma \mid c^\bullet < a^\bullet\}.
\]

Since $I\not\in\Gamma$, we have $d^\bullet < b^\bullet$ whenever $J\in \Gamma_1$. Critical points of $V$ are isolated, so there is an upper bound $s\in \gint{a}{b}$ for the endpoints of intervals in $\Gamma_1$. Consider the multiset $\Gamma_3=\{J\in \Gamma\,|\,(\pi_J)_s\neq 0\}$. This is non-empty as $\phi_s\neq 0$ by hypothesis, and $\Gamma_3\subseteq \Gamma_2$ since $(\pi_J)_s=0$ for all $J\in \Gamma_1$ by choice of $s$. Using properties of direct summands, we define 
\begin{equation*}
    y\coloneqq \beta_s(\phi_s(1_\kk)) = \sum_{J \in \Gamma} (\pi_{J})_s(1_\kk) = \sum_{J \in \Gamma_3} (\pi_{J})_s(1_\kk) \in \bigoplus_{L\in\B(V)}(\kk_L)_s.
\end{equation*}

Similarly, since critical points of $V$ are isolated, there is an upper bound $r\in(-\infty,a^\bullet\rangle$ for startpoints in $\Gamma_2$, and thus in $\Gamma_3$. In particular, $r\in J$ for all $J\in \Gamma_3$, so we can define an element
\begin{equation*}
    x\coloneqq \sum_{J\in \Gamma_3} (\pi_{J})_{r
    }(1_\kk) \in \bigoplus_{L\in\B(V)}(\kk_L)_r.
\end{equation*}
The structure map (identity) sends $x$ to $y$. But the intersection $\im V_{(r,s)}\cap \im(\phi_s)$ is non-zero, as it contains $\beta^{-1}(y)$. This contradicts our hypothesis.

We conclude that $I\in\Gamma$, and so $\pi_I\neq 0$, and the commuting diagram
\begin{equation*}
\begin{tikzcd}
    \kk_{I}\arrow[r,"\phi", hook]\arrow[d,"\pi_I",shift left] & V\arrow[d,"\cong","\beta" swap] \\
    \kk_{I}\arrow[u,"c", shift left] & 
    \kk_{I} \oplus \bigl(\textstyle\bigoplus_{L \in \B(V)\setminus I} \kk_L \bigr), \arrow[l,twoheadrightarrow] 
\end{tikzcd}
\end{equation*}
where $c_t(1_\kk)=1/(\pi_{It}(1_\kk))$ for $t\in\gint{a}{b}$, gives the isomorphism.
\end{proof}
Combining \cref{lemma: properties-delta-epsilon,lemma: phi_splits}, assuming $\alpha_r=0$ for all $r\in(-\infty,a^\bullet\rangle$, we get the commuting diagram
\begin{equation} \label{eq: split-diagram}
    \begin{tikzcd}[
    /tikz/column 3/.style={column sep=0.1cm}, 
    /tikz/column 4/.style={column sep=0.5cm}, 
    /tikz/column 5/.style={column sep=0.5cm},
    /tikz/column 6/.style={column sep=0.1cm}
    ]
        V  \ar[d, "\alpha"]  
        & 
        & 
        \Big( \arrow[ll,"\widetilde{\delta}" swap,"\cong"]
        & 
        \kk_{\gint{a}{b}} \ar[d, "1"] 
        & 
        \oplus  
        & 
        V' \ar[d,"\alpha'"]  \ar[lld, "\gamma"]
        & 
        \Big) 
        \\
        U   
        &  
        & 
        \Big( \arrow[ll,"\widetilde{\varepsilon}" swap,"\cong"] 
        & 
        \kk_{\gint{c}{d}}
        & 
        \oplus          
        & 
        U' 
        & 
        \Big),
    \end{tikzcd}
\end{equation}
where $\widetilde{\delta}|_{\kk_{\gint{a}{b}}}=\delta$ and $\widetilde{\varepsilon}|_{\kk_{\gint{c}{d}}}=\varepsilon$.

\begin{remark}\label{remark: gamma_zero_outside_intersection}
    If $\alpha_r=0$ for all $r\in(-\infty,a^\bullet\rangle$, then $\gamma_r=0$ for $r\not\in\gint{a}{d}$. Indeed, if $r\not\in {\langle c^\bullet,d^\bullet\rangle}$, then $(\kk_{\langle c^\bullet,d^\bullet\rangle})_r=0$ and so $\gamma_r=0$, and if $r\in \gint{c}{a}$, then Diagram \eqref{eq: split-diagram} reduces to $\widetilde{\varepsilon}_r \circ (\gamma_r+\alpha'_r)=\alpha_r \circ \widetilde{\delta}_r=0$, and so $\gamma_r = 0$.
   
\end{remark}

\begin{lemma}\label{lemma: gamma_zero_on_strictly_nested_bars}
    Let $\alpha:V\to U$ be a ladder module with isolated critical points such that $S(\alpha)\cong \kk_{\langle c^\bullet,\infty)\times \gint{a}{d}\times\langle b^\bullet,\infty)}\oplus A$ and $\alpha_r=0$ for all $r\in (-\infty,a^\bullet\rangle$. Using the notation of Diagram \eqref{eq: split-diagram}, 
    if an interval summand $\kk_\gint{e}{f}\subseteq V'$ has $a^\bullet \leq e^\bullet$ and $f^\bullet < b^\bullet$, and $\gamma_t|_{(\kk_\gint{e}{f})_t}\neq0$ implies $\alpha'_t|_{(\kk_\gint{e}{f})_t}=0$, then $\gamma|_{\kk_{\gint{e}{f}}}=0$.
   
\end{lemma}
\begin{proof}
    Write $R= {\langle c^\bullet,\infty)\times \gint{a}{d}\times\langle b^\bullet,\infty)}$. Fix $t\in \gint{e}{f}$ and let $\iota_t\coloneqq1_\kk\in (\kk_\gint{e}{f})_t\subseteq V'_t$. We assume by contradiction that $\zeta\coloneqq \gamma_t(\iota_t)\neq 0$ in $(\kk_\gint{c}{d})_t$. Since $\zeta$ is non-zero, we must have $t\in\gint{e}{d}\subseteq\gint{a}{d}$, and
    \begin{equation*}
        \widetilde{\varepsilon}_t(\gamma_t(\iota_t))=\varepsilon_t(\zeta) = \zeta u_t = \Psi_{(r,t,s)}(\zeta)\in \kk_R\subseteq S(\alpha)_{(r,t,s)}.
    \end{equation*}
    
    By hypothesis, we have $\alpha'_t(\iota_t)=0$, so $\zeta u_t = \widetilde{\varepsilon}_t(\gamma_t(\iota_t))= \alpha_t(\delta(\iota_t))$ by Diagram \eqref{eq: split-diagram}. By hypothesis, we can pick $s\in \gint{f}{b}$. Since $\iota_t\in\kk_{\gint{e}{f}}\subseteq V'$ we have $V'_{t,s}(\iota_t)=0$ and so $\zeta u_t=\alpha_t(\widetilde{\delta}_t(\iota_t))\in \alpha_t(\ker V_{t,s})=S(\alpha)_{(t,t,s)}$ where $(t,t,s)\not\in R$. However, now the structure map $S(\alpha)_{(t,t,s)}\hookrightarrow S(\alpha)_{(t,t,s')}$ for $s'\in\langle b^\bullet,\infty)$ includes the non-zero element $\zeta u_t\neq 0$ from outside the rectangle into the rectangle $\kk_R$, contradicting $\kk_R$ being a direct summand. \qedhere

\end{proof}

For now, we continue with the case where $S(\alpha)$ is a single rectangle and $V$ has no infinite bars.

\begin{remark}\label{remark: rectangle_simplifications} In the case when $S(\alpha)\cong \kk_R$ is a rectangle, $\alpha$ has isolated critical points and $V$ has bars of finite length, the following properties hold:
\begin{enumerate}
    \item If $\kk_\gint{e}{f}\subseteq V'$, then  $\gamma_t|_{(\kk_\gint{e}{f})_t}\neq0$ implies $\alpha'_t|_{(\kk_\gint{e}{f})_t}=0$: Let $\iota_t=1_\kk\in (\kk_\gint{e}{f})_t$ such that $\gamma_t(\iota_t)\neq 0$. Since $V$ has bars of finite length, there exists $s\in\langle f^\bullet,\infty)$ such that $\iota_t\in\ker V_{t,s}$. Assume by contradiction that $\alpha'_t(\iota_t)\neq 0$. In particular, the element $\iota_t$ hits an interval $\kk_\gint{g}{h}\subseteq U$ in addition to $\kk_\gint{c}{d}$. But now, 
    \begin{equation*}
        \dim (S(\alpha)_{(t,t,s)}) = \dim (\alpha_t(\ker V_{t,s}))\geq 2
    \end{equation*}
    contradicting that $S(\alpha)$ is a rectangle with dimension at most $1$.
    \item If $r\not\in\gint{a}{b}$, then $\alpha_r=0$. In particular $\alpha_r=0$ for all $r\in(-\infty,a^\bullet\rangle$: Assume $0\neq x =\alpha_r(y)$ for some $r\not\in\gint{a}{b}$. Since $V$ has finite bars, there exists $s\in\langle b^\bullet,\infty)$ such that $V_{r,s}(y)=0$, so $x\in \im U_{r,r}\cap \alpha_r(\ker V_{r,s})=S(\alpha)_{(r,r,s)}$, but $(r,r,s)\not\in R$, giving a contradiction.
\end{enumerate}
\end{remark}

\begin{lemma}\label{lemma: gamma_phi_lemma}
    Let $\alpha:V\to U$ be a ladder module with isolated critical points, where $V$ has bars of finite length. If $S(\alpha)\cong \kk_{\langle c^\bullet,\infty)\times \gint{a}{d}\times\langle b^\bullet,\infty)}$, then the map $\gamma:V'\to\kk_{\langle c^\bullet,d^\bullet\rangle}$ from Diagram \eqref{eq: split-diagram} factors through $\kk_{\gint{a}{b}}\xrightarrow{1}\kk_{\gint{c}{d}}$.
\end{lemma}
\begin{proof} 
    By \cref{remark: gamma_zero_outside_intersection}, the map $\gamma_t$ is only non-zero for $t\in\gint{a}{d}$. We define the map $\widetilde{\gamma}:V'\to \kk_\gint{a}{b}$ by
    \begin{equation*}
        \widetilde{\gamma}_t(x)=
        \begin{cases}
            0 & \textrm{ for }t\not\in\gint{a}{b}\\
            \gamma_t(x) & \textrm{ for } t\in\gint{a}{d}\\
            \gamma_\wde(x) & \textrm{ for }t\in\gint{d}{b}
        \end{cases}
    \end{equation*}
    where $\wde$ is such that $\alpha$ is constant on $[\wde,d\rangle$, using isolated critical points. We need to show that $\widetilde{\gamma}$ is a persistence map, so let $t\leq s$ in $\gint{a}{b}$. If both $t$ and $s$ are in $\gint{a}{d}$, then naturality follows from $\gamma$, so let $s\in\gint{d}{b}$. By construction, $V'$ is the direct sum $\bigoplus_{L\in\B(V')}\kk_L$. 
    Let $\kk_\gint{e}{f}\subseteq V'$ be such that $\gamma|_{\kk_\gint{e}{f}}\neq 0$. From \cref{remark: rectangle_simplifications,lemma: gamma_zero_on_strictly_nested_bars} we get $b^\bullet \leq f^\bullet$, noting that the case $e^\bullet < a^\bullet$ is impossible as $\alpha_r=0$ for all $r\in(-\infty,a^\bullet\rangle$. In particular, we have $s\in\gint{e}{f}$. Let $x=\sum_{L\in\B(V')}c_L$ for $c_L\in(\kk_L)_t$. If $t\in L$, then $V'_{t,s}(c_L)=c_L$ as structure maps are identities along bars, and if $t\not\in L$ then $c_L=0=V'_{t,s}(c_L)$. Thus, $V'_{t,s}(x)=x$ using linearity, and
\begin{equation*}
        \widetilde{\gamma}_s(V'_{t,s}(x))=\widetilde{\gamma}_s(x)=\gamma_{\wde}(x)=(\kk_\gint{c}{d})_{q,\wde}(\gamma_q(x))=\gamma_q(x)=\widetilde{\gamma}_t(x),
    \end{equation*}
    where $q=t$ when $t\in\langle a^\bullet,\wde)$ and $q=\wde$ when $t\in [\wde,s]$.
\end{proof}

Write $I=\langle a^\bullet,b^\bullet\rangle$ and define persistence maps $\Omega^\pm:\kk_I\oplus V'\to \kk_I\oplus V'$ by mapping points $(x,y)$ in filtration value $t$ to $(x\pm\widetilde{\gamma}_t(y),y)$. Clearly $\Omega^+$ and $\Omega^-$ are inverses of each other, and thus isomorphisms. As matrices, they are
\begin{equation*}
    \Omega^\pm = 
    \begin{bmatrix}
    \Id_{\kk_I} & \pm\widetilde{\gamma} \\
    0 & \Id_{V'}
    \end{bmatrix},
\end{equation*}
and, writing $\phi:\kk_\gint{a}{b}\xrightarrow{1}\kk_\gint{c}{d}$, we have 
\begin{equation*}
    \begin{bmatrix}
    \phi & \gamma \\
    0 & \alpha'
    \end{bmatrix}_t
    \begin{bmatrix}
    \Id_{\kk_\gint{a}{b}} & -\widetilde{\gamma} \\
    0 & \Id_{V'}
    \end{bmatrix}_t
    =
    \begin{bmatrix}
    \phi_t & \gamma_t-\phi_t\circ \widetilde{\gamma}_t \\
    0 & \alpha'_t
    \end{bmatrix}.
\end{equation*}
From \cref{remark: gamma_zero_outside_intersection}, we have $\gamma_t=0$ for $t\notin \gint{a}{d}$, where also $\phi_t=0$. If $t\in\gint{a}{b}$, then $\phi_t\circ\widetilde{\gamma}_t=\Id_\kk \circ \gamma_t=\gamma_t$. So, $\gamma_t-\phi_t\circ\widetilde{\gamma}_t=0$ for all $t$, and the following diagram commutes:
\begin{equation}\label{eq: diagram2fortheproof}
\begin{tikzcd}[column sep=0.5em]
    \kk_{\langle a^\bullet,b^\bullet\rangle} \arrow[d,"1"] & 
    \oplus           & 
    V' \arrow[d,"\alpha'"] \arrow[rrrr,"\Omega^-","\cong" swap]   &\textrm{  } &\textrm{  } &\textrm{  } & 
    \kk_{\langle a^\bullet,b^\bullet\rangle} \ar[d,"1"] & 
    \oplus           & 
    V' \ar[d,"\alpha'"]  \ar[lld, "\gamma" swap]
    \\
    \kk_{\langle c^\bullet,d^\bullet\rangle} & 
    \oplus          & 
    U'\arrow[rrrr,equal,"\Id"] & & & &
    \kk_{\langle c^\bullet,d^\bullet\rangle}& 
    \oplus          & 
    U'.
    \end{tikzcd}
\end{equation}

\converseresult*
\begin{proof}
    By \cref{prop: shape-rectangles}, we assume $R= {\langle c^\bullet,\infty)\times \gint{a}{d}\times\langle b^\bullet,\infty)}$ where $c^\bullet\leq a^\bullet\leq d^\bullet \leq b^\bullet$. Combining Diagrams \eqref{eq: split-diagram} and \eqref{eq: diagram2fortheproof}, we get $\alpha\cong \phi\oplus\alpha'$, where $\phi:\kk_\gint{a}{b}\to\kk_\gint{c}{d}$ is a non-trivial summand of $\alpha$. 

    We finish by showing that $\alpha'$ is trivial. Assume by contradiction that there is a $y\in V'_t$ such that $0\neq\alpha'_t(y)=x$. Since, $V$ (and thus $V')$ has bars of finite length, there exists $s\geq t$ such that $V'_{t,s}(y)=0$. In particular, we have a non-empty element $x$ in $\alpha'_t(\ker V'_{t,s})= S(\alpha')_{(t,t,s)}$. By additivity of $S$ (\cref{prop: S_additive}), we get $S(\alpha)\cong S(\phi)\oplus S(\alpha')$, so $S(\alpha)$ is a direct sum of two non-empty modules, contradicting the hypothesis that it is a rectangle. 
\end{proof}

\section*{Acknowledgements}
Lars M. Salbu was partially co-funded by the European Union, grant GA\#101126560; Bergen research and training
program for future AI leaders across the disciplines, LEAD AI.

Marco Delgado Garrido was supported by grant PREX2024-000254 funded by MICIU/AEI/10.13039/501100011033 and European Union (ESF+) (IMUS María de Maeztu Excellence Unit grant CEX2024-001517-M).
\end{document}